\documentclass[a4paper,11pt,twoside,fleqn]{article}
\usepackage{amsmath,amssymb,amsthm,latexsym}
\newtheorem{theorem}{Theorem}[section]
\newtheorem{proposition}[theorem]{Proposition}
\newtheorem{corollary}[theorem]{Corollary}
\newtheorem{lemma}[theorem]{Lemma}

\renewcommand{\Re}{\text{Re}}

\newcommand{\dist}{\text{\rm dist}}
\newcommand{\Id}{\text{\rm I}}
\newcommand{\proj}{\text{\rm proj}\,}

\begin{document}

\title{Complex anti-self-dual instantons and Cayley submanifolds}
\author{Simon Brendle}
\date{August 13, 2003}

\maketitle

\section{Introduction}

Let $M$ be a manifold of dimension $8$, and let $\Omega$ be a $4$-form which defines an almost 
$Spin(7)$-structure on $M$. An $\Omega$-anti-self-dual instanton is a connection $A$ on a vector 
bundle over $M$ such that the curvature $F_A$ satisfies 
\begin{equation} 
F_A + *(\Omega \wedge F_A) = 0. 
\end{equation} 
If $M$ is an almost Calabi-Yau manifold, then the $4$-form $\Omega$ can be written as 
\[\Omega = 4 \, \Re(\theta) + \frac{1}{2} \, \omega^2,\] 
where $\omega \in \Omega^{1,1}(M)$ denotes the symplectic form and $\theta \in \Omega^{0,4}(M)$ 
is the complex volume form. The complex volume form induces an anti-linear involution $*_\theta: 
\Omega^{0,2}(M) \to \Omega^{0,2}(M)$. Then the anti-self-duality equation (1) is equivalent to 
\begin{equation} 
F_A^{1,1} \cdot \omega = 0 
\end{equation} 
and 
\begin{equation} 
(1 + *_\theta) \, F_A^{0,2} = 0. 
\end{equation} 
The space of $2$-forms splits as a direct sum 
\begin{equation} 
\Lambda^2 TM = \Lambda_+^2 TM \oplus \Lambda_-^2 TM, 
\end{equation} 
where 
\begin{equation} 
\Lambda_+^2 TM = \{\varphi \in \Lambda^2 M: 3\varphi - *(\Omega \wedge \varphi) = 0\} 
\end{equation}  
and 
\begin{equation} 
\Lambda_-^2 TM = \{\varphi \in \Lambda^2 M: \varphi + *(\Omega \wedge \varphi) = 0\}. 
\end{equation} 
Note that $\Lambda_+^2 M$ is a vector space of dimension $7$ and $\Lambda_-^2(M)$ is a vector 
space of dimension $21$. Let $P_+$ and $P_-$ be the projections associated to the splitting (4). 
This implies 
\[P_+ \varphi = \frac{1}{4} \, (\varphi + *(\Omega \wedge \varphi))\] 
and 
\[P_- \varphi = \frac{1}{4} \, (3\varphi - *(\Omega \wedge \varphi)).\] 
We denote by $\Omega_+^2(M)$ the space of sections of the vector bundle $\Lambda_+^2 TM$. 
Similarly, $\Omega_-^2(M)$ is the space of sections of the vector bundle $\Lambda_-^2 TM$. \\

If $\Omega$ is closed, then the anti-self-duality equation (1) implies the Yang-Mills equation 
$D_A^* F_A = 0$. \\

The equations (1),(2) generalize the anti-self-dual equations in dimension $4$ (see e.g. \cite
{DK,Ta1}), and have been studied by various authors, including S. K. Donaldson and R. P. Thomas 
\cite{DT,Th}, L. Baulieu, H. Kanno, and I. M. Singer \cite{BKS}, J. Chen \cite{Ch}, and G. Tian 
\cite{Ti}. These submanifolds are also of considerable interest in mathematical physics. \\

G. Tian constructed a compactification of the moduli space of $\Omega$-anti-self-dual instantons 
over $M$. He proved that every sequence $A_k$ of $\Omega$-anti-self-dual instantons over $M$ has 
a subsequence, still denoted by $A_k$, such that 
\[\lim_{k \to \infty} \int_M c_2(A_k) \wedge \psi = \int_M c_2(A_\infty) \wedge \psi + \int_S 
\Theta \, \psi,\] 
where $c_2$ denotes the $4$-form representing the second Chern class of the bundle, and $\psi$ 
is a smooth $4$-form on $M$. Furthermore, $A_\infty$ is a $\Omega$-anti-self-dual instanton 
which is smooth outside a set of vanishing $\mathcal{H}^4$-measure. Furthermore, $S$ is a Cayley 
submanifold, i.e. a submanifold calibrated by the $4$-form $\Omega$. Cayley submanifolds were 
studied by R. Harvey and H. B. Lawson \cite{HL}. There is a rich class of examples. For 
instance, this class contains as limiting cases the holomorphic subvarieties and the special 
Lagrangian submanifolds of $M$. Special Lagrangian submanifolds have been studied extensively, 
see e.g. \cite{Hi}. Cayley submanifolds play a role in high-energy physics, see for example 
\cite{BBMOOY}. \\

Our aim in this paper is to construct smooth complex anti-self-dual instantons such that the 
energy density $|F_A|^2$ is concentrated near a given Cayley submanifold $S$. \\

In the first step, we construct a suitable family of approximate solutions. To this end, we 
assume that the normal bundle $NS$ can be endowed with a complex structure $J$ and a complex 
volume form $\omega$. Each approximate solution is described by a set $(v,\lambda,J,\omega)$, 
where $v$ is a section of the normal bundle of $S$, $\lambda$ is a positive function on $S$, and 
$(J,\omega)$ is a $SU(2)$-structure on $NS$. The covariant derivative of the pair $(J,\omega)$ 
can be described by a $1$-form $\theta$ with values in the Lie algebra $\Lambda_+^2 NS$. \\

The covariant derivative of the $4$-form $\Omega$ can be written in the form 
\[\nabla_X \Omega = \sum_{k=1}^8 i_{e_k} \alpha(X) \wedge i_{e_k} \Omega,\] 
where $\alpha$ is a $1$-form with values in $\Lambda_+^2 TM$. \\

We consider the elliptic complex 
\[0 \longrightarrow \Omega^0(M) \begin{array}[b]{c} d \\ \longrightarrow \end{array} \Omega^1(M) 
\begin{array}[b]{c} P_+ d \\ \longrightarrow \end{array} \Omega_+^2(M) \longrightarrow 0.\] 
The first and the second cohomology groups associated to this elliptic comples are $H^0(M)$ and 
$H^1(M)$. The third cohomology group is denoted by $H_+^2(M)$. \\

\begin{theorem}
Suppose that $H_+^2(M) = 0$. Then, for each $\varepsilon > 0$, there exists a mapping $\Xi
_\varepsilon$ whch assigns to each set of glueing data $(v,\lambda,J,\omega) \in \mathcal{C}^{2,
\gamma}(S)$ a section of the vector bundle $V \oplus W$ of class $\mathcal{C}^\gamma(S)$ such 
that the following holds. \\

(i) If $(v,\lambda,J,\omega)$ is a set of glueing data such that 
\[\|v\|_{\mathcal{C}^{1,\gamma}(S)} \leq K,\] 
\[\|\lambda\|_{\mathcal{C}^{1,\gamma}(S)} \leq K, \qquad \inf \lambda \geq 1,\] 
\[\|(J,\omega)\|_{\mathcal{C}^{1,\gamma}(S)} \leq K,\] 
then we have the estimate 
\begin{align*} 
\bigg \| &\Xi_\varepsilon(v,\lambda,J,\omega) \\ 
&- 4 \, \bigg ( \proj_V \Big ( \sum_{i,j=1}^4 (\nabla_i v_k + \alpha_{ik,l} \, v_l) \, e_i 
\otimes e_k^\perp \Big ), \\ 
&\hspace{10.1mm} \proj_W \Big ( \sum_{i,k,l=1}^4 (\lambda^{-1} \, \nabla_i \lambda \, \delta
_{kl} + \theta_{i,kl} + \alpha_{ik,l}) \, e_i \otimes e_k^\perp \otimes e_l^\perp \bigg ) \bigg 
\|_{\mathcal{C}^\gamma(S)} \leq C \, \varepsilon^{\frac{1}{32}}. 
\end{align*}

(ii) If $\Xi_\varepsilon(v,\lambda,J,\omega) = 0$, then the approximate solution $A$ 
corresponding to $(v,\lambda,J,\omega)$ can be deformed to a nearby connection $\tilde{A}$ 
satisfying $F_{\tilde{A}} + *(\Omega \wedge F_{\tilde{A}}) = 0$.
\end{theorem}

\vspace{2mm}

In Section 2, we study the mapping properties of a model operator on $\mathbb{R}^8$. \\

In Section 3, we construct a family of approximate solutions of the Yang-Mills equations. 
More precisely, given any set of glueing data $(v,\lambda,J,\omega)$ satisfying 
\[\|v\|_{\mathcal{C}^{1,\gamma}(S)} \leq K,\] 
\[\|\lambda\|_{\mathcal{C}^{1,\gamma}(S)} \leq K, \qquad \inf \lambda \geq 1,\] 
\[\|(J,\omega)\|_{\mathcal{C}^{1,\gamma}(S)} \leq K,\] 
we construct a connection $A$ such that 
\[\|F_A + *(\Omega \wedge F_A)\|_{\mathcal{C}_3^\gamma(M)} \leq C \, \varepsilon^2.\] 
Here, the weighted H\"older space $\mathcal{C}_\nu^\gamma(M)$ is defined as 
\begin{align*} 
\|u\|_{\mathcal{C}_\nu^\gamma(M)} 
&= \sup \, (\varepsilon + \dist(p,S))^\nu \, |u(p)| \\ 
&+ \sup_{\begin{smallmatrix} 4\dist(p_1,p_2) \leq \\ \varepsilon + \dist(p_1,S) + \dist(p_2,S) 
\end{smallmatrix}} \, (\varepsilon + \dist(p_1,S) + \dist(p_2,S))^{\nu+\gamma} \frac{|u(p_1) - 
u(p_2)|}{\dist(p_1,p_2)^\gamma}. 
\end{align*} 

In Section 4, we derive estimates for the linearized operator which are independent of 
$\varepsilon$. \\

In Section 5, we apply the contraction mapping principle to deform the approximate solution $A$ 
to a nearby connection $\tilde{A} = A + a$ such that 
\[(\Id - \mathbb{P}) (F_{\tilde{A}} + *(\Omega \wedge F_{\tilde{A}})) = 0,\] 
where $(\Id - \mathbb{P})$ is the fibrewise projection from $\mathcal{C}_\nu^\gamma(M)$ to the 
subspace $\mathcal{G}_\nu^\gamma(M)$. In particular, if the balancing condition 
\[\mathbb{P} (F_{\tilde{A}} + *(\Omega \wedge F_{\tilde{A}})) = 0\] 
is satisfied, then $\tilde{A}$ is an $\Omega$-anti-self-dual instanton. \\

In Section 6, we calculate the leading term in the asymptotic expansion of  
\[\mathbb{P} (F_{\tilde{A}} + *(\Omega \wedge F_{\tilde{A}})) = 0.\] 
This concludes the proof of Theorem 1.1. \\

The author is grateful to Professor Gerhard Huisken and Professor Gang Tian for discussions. \\

\section{The model problem on $\mathbb{R}^8$}

The $Spin(7)$-structure on $\mathbb{R}^8$ is given by 
\begin{align*} 
\Omega 
&= -e_1 \wedge e_2 \wedge e_1^\perp \wedge e_2^\perp - e_1 \wedge e_2 \wedge e_3^\perp \wedge 
e_4^\perp - e_3 \wedge e_4 \wedge e_1^\perp \wedge e_2^\perp \\ 
&- e_3 \wedge e_4 \wedge e_3^\perp \wedge e_4^\perp + e_1 \wedge e_3 \wedge e_2^\perp \wedge e_4
^\perp - e_1 \wedge e_3 \wedge e_1^\perp \wedge e_3^\perp \\ 
&- e_2 \wedge e_4 \wedge e_2^\perp \wedge e_4^\perp + e_2 \wedge e_4 \wedge e_1^\perp \wedge e_3
^\perp - e_1 \wedge e_4 \wedge e_2^\perp \wedge e_3^\perp \\ 
&- e_1 \wedge e_4 \wedge e_1^\perp \wedge e_4^\perp - e_2 \wedge e_3 \wedge e_2^\perp \wedge e_3
^\perp - e_2 \wedge e_3 \wedge e_1^\perp \wedge e_4^\perp \\ 
&+ e_1 \wedge e_2 \wedge e_3 \wedge e_4 + e_1^\perp \wedge e_2^\perp \wedge e_3^\perp \wedge e_4
^\perp. 
\end{align*} 
Hence, the $2$-forms 
\begin{align*} 
&e_1 \wedge e_2 + e_3 \wedge e_4 - e_1^\perp \wedge e_2^\perp - e_3^\perp \wedge e_4^\perp, \\ 
&e_1 \wedge e_3 - e_2 \wedge e_4 - e_1^\perp \wedge e_3^\perp + e_2^\perp \wedge e_4^\perp, \\ 
&e_1 \wedge e_4 + e_2 \wedge e_3 - e_1^\perp \wedge e_4^\perp - e_2^\perp \wedge e_3^\perp, \\ 
&e_1 \wedge e_1^\perp + e_2 \wedge e_2^\perp + e_3 \wedge e_3^\perp + e_4 \wedge e_4^\perp, \\ 
&e_1 \wedge e_2^\perp - e_2 \wedge e_1^\perp - e_3 \wedge e_4^\perp + e_4 \wedge e_3^\perp, \\ 
&e_1 \wedge e_3^\perp + e_2 \wedge e_4^\perp - e_3 \wedge e_1^\perp - e_4 \wedge e_2^\perp, \\ 
&e_1 \wedge e_4^\perp - e_2 \wedge e_3^\perp + e_3 \wedge e_2^\perp - e_4 \wedge e_1^\perp 
\end{align*} 
form a basis for $\Lambda_+^2 \mathbb{R}^8$. \\

Let now $E$ be a Cayley subspace of $\mathbb{R}^8$, i.e. a subspace calibrated by $\Omega$. The 
group $Spin(7)$ acts transitively on the set of Cayley subspaces and leaves the $4$-form 
$\Omega$ invariant (cf. \cite{HL}). Hence, we may assume without loss of generality that $E$ is 
spanned by $\{e_i: 1 \leq i \leq 4\}$ and $E^\perp$ is spanned by $\{e_i^\perp: 1 \leq i \leq 4
\}$. \\

We define two vector spaces $V \subset E \otimes E^\perp$ and $W \subset E \otimes E^\perp 
\otimes E^\perp$ over the submanifold $S$. The following elements form a basis for $V$: 
\[e_1 \otimes e_1^\perp + e_2 \otimes e_2^\perp + e_3 \otimes e_3^\perp + e_4 \otimes e_4^\perp,
\] 
\[e_1 \otimes e_2^\perp - e_2 \otimes e_1^\perp - e_3 \otimes e_4^\perp + e_4 \otimes e_3^\perp,
\] 
\[e_1 \otimes e_3^\perp + e_2 \otimes e_4^\perp - e_3 \otimes e_1^\perp - e_4 \otimes e_2^\perp,
\] 
\[e_1 \otimes e_4^\perp - e_2 \otimes e_3^\perp + e_3 \otimes e_2^\perp - e_4 \otimes e_1^\perp.
\] 
The following elements form a basis for $W$: 
\begin{align*} 
&e_1 \otimes (e_1^\perp \otimes e_1^\perp + e_2^\perp \otimes e_2^\perp + e_3^\perp \otimes e
_3^\perp + e_4^\perp \otimes e_4^\perp) \\ 
&+ e_2 \otimes (e_2^\perp \otimes e_1^\perp - e_1^\perp \otimes e_2^\perp + e_4^\perp \otimes e
_3^\perp - e_3^\perp \otimes e_4^\perp) \\ 
&+ e_3 \otimes (e_3^\perp \otimes e_1^\perp - e_4^\perp \otimes e_2^\perp - e_1^\perp \otimes e
_3^\perp + e_2^\perp \otimes e_4^\perp) \\ 
&+ e_4 \otimes (e_4^\perp \otimes e_1^\perp + e_3^\perp \otimes e_2^\perp - e_2^\perp \otimes e
_3^\perp - e_1^\perp \otimes e_4^\perp), 
\end{align*} 
\begin{align*} 
&e_1 \otimes (e_2^\perp \otimes e_1^\perp - e_1^\perp \otimes e_2^\perp + e_4^\perp \otimes e
_3^\perp - e_3^\perp \otimes e_4^\perp) \\ 
&- e_2 \otimes (e_1^\perp \otimes e_1^\perp + e_2^\perp \otimes e_2^\perp + e_3^\perp \otimes e
_3^\perp + e_4^\perp \otimes e_4^\perp) \\ 
&- e_3 \otimes (e_4^\perp \otimes e_1^\perp + e_3^\perp \otimes e_2^\perp - e_2^\perp \otimes e
_3^\perp - e_1^\perp \otimes e_4^\perp) \\  
&+ e_4 \otimes (e_3^\perp \otimes e_1^\perp - e_4^\perp \otimes e_2^\perp - e_1^\perp \otimes e
_3^\perp + e_2^\perp \otimes e_4^\perp), 
\end{align*}
\begin{align*} 
&e_1 \otimes (e_3^\perp \otimes e_1^\perp - e_4^\perp \otimes e_2^\perp - e_1^\perp \otimes e
_3^\perp + e_2^\perp \otimes e_4^\perp) \\ 
&+ e_2 \otimes (e_4^\perp \otimes e_1^\perp + e_3^\perp \otimes e_2^\perp - e_2^\perp \otimes e
_3^\perp - e_1^\perp \otimes e_4^\perp) \\ 
&- e_3 \otimes (e_1^\perp \otimes e_1^\perp + e_2^\perp \otimes e_2^\perp + e_3^\perp \otimes e
_3^\perp + e_4^\perp \otimes e_4^\perp) \\ 
&- e_4 \otimes (e_2^\perp \otimes e_1^\perp - e_1^\perp \otimes e_2^\perp + e_4^\perp \otimes e
_3^\perp - e_3^\perp \otimes e_4^\perp), 
\end{align*} 
\begin{align*}
&e_1 \otimes (e_4^\perp \otimes e_1^\perp + e_3^\perp \otimes e_2^\perp - e_2^\perp \otimes e
_3^\perp - e_1^\perp \otimes e_4^\perp) \\ 
&- e_2 \otimes (e_3^\perp \otimes e_1^\perp - e_4^\perp \otimes e_2^\perp - e_1^\perp \otimes e
_3^\perp + e_2^\perp \otimes e_4^\perp) \\ 
&+ e_3 \otimes (e_2^\perp \otimes e_1^\perp - e_1^\perp \otimes e_2^\perp + e_4^\perp \otimes e
_3^\perp - e_3^\perp \otimes e_4^\perp) \\ 
&- e_4 \otimes (e_1^\perp \otimes e_1^\perp + e_2^\perp \otimes e_2^\perp + e_3^\perp \otimes e
_3^\perp + e_4^\perp \otimes e_4^\perp). 
\end{align*} 

Let $B$ be a connection which is invariant under translations along $E$ and agrees with the 
basic instanton along $E^\perp$. More precisely, we define 
\begin{align*} 
B(e_1^\perp) &= \frac{-y_2 \, \mathfrak{i} - y_3 \, \mathfrak{j} - y_4 \, \mathfrak
{k}}{\varepsilon^2 + |y|^2} \\ 
B(e_2^\perp) &= \frac{y_1 \, \mathfrak{i} - y_4 \, \mathfrak{j} + y_3 \, \mathfrak
{k}}{\varepsilon^2 + |y|^2} \\ 
B(e_3^\perp) &= \frac{y_4 \, \mathfrak{i} + y_1 \, \mathfrak{j} - y_2 \, \mathfrak
{k}}{\varepsilon^2 + |y|^2} \\ 
B(e_4^\perp) &= \frac{-y_3 \, \mathfrak{i} + y_2 \, \mathfrak{j} + y_1 \, \mathfrak
{k}}{\varepsilon^2 + |y|^2}, 
\end{align*} 
where 
\begin{align*} 
&\mathfrak{i}(e_1^\perp) = -e_2^\perp, \quad \mathfrak{i}(e_2^\perp) = e_1^\perp, \quad 
\mathfrak{i}(e_3^\perp) = e_4^\perp, \quad \mathfrak{i}(e_4^\perp) = -e_3^\perp, \\ 
&\mathfrak{j}(e_1^\perp) = -e_3^\perp, \quad \mathfrak{j}(e_2^\perp) = -e_4^\perp, \quad 
\mathfrak{j}(e_3^\perp) = e_1^\perp, \quad \mathfrak{j}(e_4^\perp) = e_2^\perp, \\ 
&\mathfrak{k}(e_1^\perp) = -e_4^\perp, \quad \mathfrak{k}(e_2^\perp) = e_3^\perp, \quad 
\mathfrak{k}(e_3^\perp) = -e_2^\perp, \quad \mathfrak{k}(e_4^\perp) = e_1^\perp. 
\end{align*} 
Furthermore, $B(e_i) = 0$ for $1 \leq i \leq 4$. \\

The linearized operator $\mathbb{L}_B: \Omega^1(\mathbb{R}^8) \to \Omega_+^2(\mathbb{R}^8)$ is 
given by 
\[\mathbb{L}_B a = 2 \, P_+ D_B a.\] 
The adjoint operator $\mathbb{L}_B^*: \Omega_+^2(\mathbb{R}^8) \to \Omega^1(\mathbb{R}^8)$ is 
given by 
\[\mathbb{L}_B^* \varphi = 2 \, D_B^* \varphi.\] 

We define the weighted H\"older space $\mathcal{C}_\nu^\gamma(\mathbb{R}^8)$ by 
\begin{align*} 
\|u\|_{\mathcal{C}_\nu^\gamma(\mathbb{R}^8)} 
&= \sup \, (\varepsilon + |y|)^\nu \, |u(x,y)| \\ 
&+ \sup_{\begin{smallmatrix} 4 (|x_1 - x_2| + |y_1 - y_2|) \leq \\ \varepsilon + |y_1| + |y_2| 
\end{smallmatrix}} \, (\varepsilon + |y_1| + |y_2|)^{\nu+\gamma} \frac{|u(x_1,y_1) - u(x_2,y
_2)|}{(|x_1 - x_2| + |y_1 - y_2|)^\gamma}. 
\end{align*} 
More generally, we define 
\[\|u\|_{\mathcal{C}_\nu^{k,\gamma}(\mathbb{R}^8)} = \sum_{l=0}^k \|\nabla^l u\|_{\mathcal{C}
_{\nu+l}^\gamma(\mathbb{R}^8)}.\] 

Let $\mathcal{G}_\nu^{k,\gamma}(\mathbb{R}^8)$ be the set of all $\varphi \in \Omega_+^2(\mathbb
{R}^8)$ such that $\varphi \in \mathcal{C}_\nu^{k,\gamma}(\mathbb{R}^8)$, and 
\[\int_{x + E^\perp} \sum_{i,j=1}^4 \big ( \varepsilon \, s_{ik} + t_{ikl} \, y_l \big ) \, 
\langle \varphi(e_i,e_j^\perp),F_B(e_k^\perp,e_j^\perp) \rangle = 0\] 
for all $x \in E$, $s \in V$, and $t \in W$. \\

We first derive a Weitzenb\"ock formula for the operator $\mathbb{L}_B \mathbb{L}_B^*: \Omega_+
^2(\mathbb{R}^8) \to \Omega_+^2(\mathbb{R}^8)$. We shall need two algebraic facts which can be 
verified by direct calculation. For simplicity, let $e_5 = e_1^\perp$, $e_6 = e_2^\perp$, $e_7 = 
e_3^\perp$, $e_8 = e_4^\perp$. \\

\begin{lemma}
For every $\varphi \in \Lambda_+^2 \mathbb{R}^8$, we have 
\[2 \, P_+ \big ( e_k \wedge (i_{e_l} \, \varphi) + e_l \wedge (i_{e_k} \, \varphi) \big ) = 
\delta_{kl} \, \varphi.\]
\end{lemma}

\vspace{2mm}

\begin{lemma}
For every $\varphi \in \Lambda_+^2 \mathbb{R}^8$, we have 
\[\sum_{k,l=1}^8 e_k \wedge [F_B(e_k,e_l),i_{e_l} \varphi] \in \Lambda_+^2 \mathbb{R}^8.\] 
\end{lemma}

\vspace{2mm}

\begin{proposition}
The operator $\mathbb{L}_B \mathbb{L}_B^*$ satisfies the Weitzenb\"ock formula 
\[\mathbb{L}_B \mathbb{L}_B^* \varphi = \nabla_B^* \nabla_B \varphi - 2 \sum_{k,l=1}^4 e_k^\perp 
\wedge [F_B(e_k^\perp,e_l^\perp),i_{e_l^\perp} \varphi]\] 
for every $\varphi \in \Omega_+^2(\mathbb{R}^8)$.
\end{proposition}

\textbf{Proof.} 
For every $\varphi \in \Omega_+^2(\mathbb{R}^8)$, we obtain 
\begin{align*} 
4 \, D_B D_B^* \varphi 
&= -4 \sum_{k,l=1}^8 e_k \wedge (i_{e_l} \, D_{B,e_k} D_{B,e_l} \varphi) \\ 
&= -2 \sum_{k,l=1}^8 \big ( e_k \wedge (i_{e_l} \, D_{B,e_k} D_{B,e_l} \varphi) + e_l \wedge (i
_{e_k} \, D_{B,e_k} D_{B,e_l} \varphi) \big ) \\ 
&- 2 \sum_{k,l=1}^8 e_k \wedge \big ( i_{e_l} \, (D_{B,e_k} D_{B,e_l} \varphi - D_{B,e_l} D_{B,e
_k} \varphi) \big ) \\ 
&= -2 \sum_{k,l=1}^8 \big ( e_k \wedge (i_{e_l} \, D_{B,e_k} D_{B,e_l} \varphi) + e_l \wedge (i
_{e_k} \, D_{B,e_k} D_{B,e_l} \varphi) \big ) \\ 
&- 2 \sum_{k,l=1}^8 e_k \wedge [F_B(e_k,e_l),i_{e_l} \varphi]. 
\end{align*} 
Since $\varphi \in \Omega_+^2(\mathbb{R}^8)$ and $\Omega$ is parallel, it follows that $D_{B,e
_k} D_{B,e_l} \varphi \in \Omega_+^2(\mathbb{R}^8)$. Using Lemma 2.1, we obtain 
\[2 \, P_+ \big ( e_k \wedge (i_{e_l} \, D_{B,e_k} D_{B,e_l} \varphi) + e_l \wedge (i_{e_k} \, 
D_{B,e_k} D_{B,e_l} \varphi) \big ) = \delta_{kl} \, D_{B,e_k} D_{B,e_l} \varphi.\] 
From this it follows that 
\[4 \, P_+ D_B D_B^* \varphi = \nabla_B^* \nabla_B \varphi - 2 \sum_{k,l=1}^8 P_+(e_k \wedge [F
_B(e_k,e_l),i_{e_l} \varphi]).\] 
Moreover, Lemma 2.2 implies that 
\[\sum_{k,l=1}^8 e_k \wedge [F_B(e_k,e_l),i_{e_l} \varphi] \in \Omega_+^2(\mathbb{R}^8).\] 
Thus, we conclude that 
\[4 \, P_+ D_B D_B^* \varphi = \nabla_B^* \nabla_B \varphi - 2 \sum_{k,l=1}^8 e_k \wedge [F_B(e
_k,e_l),i_{e_l} \varphi].\] 
This proves the assertion. \\

\begin{proposition}
Suppose that $\psi \in \mathcal{G}_{3+\nu}^\gamma(\mathbb{R}^8)$ has compact support. Then there 
exists some $\varphi \in \mathcal{G}_{1+\nu}^{2,\gamma}(\mathbb{R}^8)$ such that 
\[\|\varphi\|_{\mathcal{C}_{1+\nu}^{2,\gamma}(\mathbb{R}^8)} \leq C \, \|\psi\|_{\mathcal{C}_{3+
\nu}^\gamma(\mathbb{R}^8)}\] 
and 
\[\mathbb{L}_B \mathbb{L}_B^* \varphi = \psi.\] 
\end{proposition}

\textbf{Proof.} 
Since $\psi \in \Omega_+^2(\mathbb{R}^8)$, we may write 
\begin{align*} 
\psi 
&= (e_1 \wedge e_2 + e_3 \wedge e_4 - e_1^\perp \wedge e_2^\perp - e_3^\perp \wedge e_4^\perp) 
\otimes g_2 \\ 
&+ (e_1 \wedge e_3 - e_2 \wedge e_4 - e_1^\perp \wedge e_3^\perp + e_2^\perp \wedge e_4^\perp) 
\otimes g_3 \\ 
&+ (e_1 \wedge e_4 + e_2 \wedge e_3 - e_1^\perp \wedge e_4^\perp - e_2^\perp \wedge e_3^\perp) 
\otimes g_4 \\ 
&+ (e_1 \wedge e_1^\perp + e_2 \wedge e_2^\perp + e_3 \wedge e_3^\perp + e_4 \wedge e_4^\perp) 
\otimes g_1^\perp \\ 
&+ (e_1 \wedge e_2^\perp - e_2 \wedge e_1^\perp - e_3 \wedge e_4^\perp + e_4 \wedge e_3^\perp) 
\otimes g_2^\perp \\ 
&+ (e_1 \wedge e_3^\perp + e_2 \wedge e_4^\perp - e_3 \wedge e_1^\perp - e_4 \wedge e_2^\perp) 
\otimes g_3^\perp \\ 
&+ (e_1 \wedge e_4^\perp - e_2 \wedge e_3^\perp + e_3 \wedge e_2^\perp - e_4 \wedge e_1^\perp) 
\otimes g_4^\perp, 
\end{align*} 
where $g_j,g_j^\perp \in \mathcal{C}_{3+\nu}^\gamma(\mathbb{R}^8)$. Furthermore, since $\psi \in 
\mathcal{G}_{3+\nu}^\gamma(\mathbb{R}^8)$, we deduce that 
\[\int_{x + E^\perp} \sum_{j=1}^4 \langle g_j^\perp,F_B(X,e_j^\perp) \rangle = 0\] 
for all $x \in E$ and all vector fields of the form 
\[X = \varepsilon \, w_k \, e_k^\perp + \mu \, y_k \, e_k^\perp + r_{kl} \, y_l \, e_k^\perp.\] 
Using Corollary 3.6 in \cite{Br}, we can find $f_j,f_j^\perp \in \mathcal{C}_{1+\nu}^{2,\gamma}
(\mathbb{R}^8)$ such that 
\[\nabla_B^* \nabla_B f_j = g_j\] 
and 
\[\nabla_B^* \nabla_B f_j^\perp - 2 \sum_{k=1}^4 [F_B(e_j^\perp,e_k^\perp),f_k^\perp] = g_j
^\perp\] 
and 
\[\int_{x + E^\perp} \sum_{j=1}^4 \langle f_j^\perp,F_B(X,e_j^\perp) \rangle = 0\] 
for all $x \in E$ and all vector fields of the form 
\[X = \varepsilon \, w_k \, e_k^\perp + \mu \, y_k \, e_k^\perp + r_{kl} \, y_l \, e_k^\perp.\] 
We now define 
\begin{align*} 
\varphi 
&= (e_1 \wedge e_2 + e_3 \wedge e_4 - e_1^\perp \wedge e_2^\perp - e_3^\perp \wedge e_4^\perp) 
\otimes f_2 \\ 
&+ (e_1 \wedge e_3 - e_2 \wedge e_4 - e_1^\perp \wedge e_3^\perp + e_2^\perp \wedge e_4^\perp) 
\otimes f_3 \\ 
&+ (e_1 \wedge e_4 + e_2 \wedge e_3 - e_1^\perp \wedge e_4^\perp - e_2^\perp \wedge e_3^\perp) 
\otimes f_4 \\ 
&+ (e_1 \wedge e_1^\perp + e_2 \wedge e_2^\perp + e_3 \wedge e_3^\perp + e_4 \wedge e_4^\perp) 
\otimes f_1^\perp \\ 
&+ (e_1 \wedge e_2^\perp - e_2 \wedge e_1^\perp - e_3 \wedge e_4^\perp + e_4 \wedge e_3^\perp) 
\otimes f_2^\perp \\ 
&+ (e_1 \wedge e_3^\perp + e_2 \wedge e_4^\perp - e_3 \wedge e_1^\perp - e_4 \wedge e_2^\perp) 
\otimes f_3^\perp \\ 
&+ (e_1 \wedge e_4^\perp - e_2 \wedge e_3^\perp + e_3 \wedge e_2^\perp - e_4 \wedge e_1^\perp) 
\otimes f_4^\perp. 
\end{align*} 
Then $\varphi \in \mathcal{G}_{1+\nu}^{2,\gamma}(\mathbb{R}^8)$, and 
\[\mathbb{L}_B \mathbb{L}_B^* \varphi = \psi\] 
by Proposition 2.3. This proves the assertion. \\

\begin{corollary}
Suppose that $\psi \in \mathcal{G}_{3+\nu}^\gamma(\mathbb{R}^8)$ has compact support. Then there 
exists a $1$-form $a \in \mathcal{C}_{2+\nu}^{1,\gamma}(\mathbb{R}^8)$ such that 
\[\|a\|_{\mathcal{C}_{2+\nu}^{1,\gamma}(\mathbb{R}^8)} \leq C \, \|\psi\|_{\mathcal{C}_{3+\nu}
^\gamma(\mathbb{R}^8)}\] 
and 
\[\mathbb{L}_B a = \psi.\] 
\end{corollary}

\textbf{Proof.} 
By Proposition 2.4, there exists some $\varphi \in \mathcal{G}_{1+\nu}^{2,\gamma}(\mathbb{R}^8)$ 
such that 
\[\|\varphi\|_{\mathcal{C}_{1+\nu}^{2,\gamma}(\mathbb{R}^8)} \leq C \, \|\psi\|_{\mathcal{C}_{3+
\nu}^\gamma(\mathbb{R}^8)}\] 
and 
\[\mathbb{L}_B \mathbb{L}_B^* \varphi = \psi.\] 
Let $a = \mathbb{L}_B^* \varphi$. Then $a \in \mathcal{C}_{2+\nu}^{1,\gamma}(\mathbb{R}^8)$ 
satisfies 
\[\|a\|_{\mathcal{C}_{2+\nu}^{1,\gamma}(\mathbb{R}^8)} \leq C \, \|\psi\|_{\mathcal{C}_{3+\nu}
^\gamma(\mathbb{R}^8)}\] 
and 
\[\mathbb{L}_B a = \psi.\] 
This proves the assertion. \\

\begin{proposition}
Let $0 < \nu < 1$. Suppose that $\psi \in \mathcal{G}_{3+\nu}^\gamma(\mathbb{R}^8)$ is supported 
in the set $\{(x,y) \in \mathbb{R}^8: |x| \leq \delta, \, |y| \leq 2\delta^4\}$. Then there 
exists a $1$-form $a \in \mathcal{C}_{2+\nu}^{1,\gamma}(\mathbb{R}^8)$ such that $a$ is 
supported in $\{(x,y) \in \mathbb{R}^8: |x| \leq 2\delta, \, |y| \leq 2\delta^2\}$, 
\[\|a\|_{\mathcal{C}_{2+\nu}^{1,\gamma}(\mathbb{R}^8)} \leq C \, \|\psi\|_{\mathcal{C}_{3+\nu}
^\gamma(\mathbb{R}^8)}\] 
and 
\[\|\mathbb{L}_B a - \psi\|_{\mathcal{C}_{3+\nu}^\gamma(\{(x,y) \in \mathbb{R}^8: |y| \leq 
2\delta^4\})} \leq C \, \delta \, \|\psi\|_{\mathcal{C}_{3+\nu}^\gamma(\mathbb{R}^8)},\] 
and 
\[\|\mathbb{L}_B a - \psi\|_{\mathcal{C}_{3+\nu}^\gamma(\mathbb{R}^8)} \leq C \, |\log \delta|
^{-1} \, \|\psi\|_{\mathcal{C}_{3+\nu}^\gamma(\mathbb{R}^8)}.\] 
\end{proposition}

\textbf{Proof.} 
By Corollary 2.5, there exists a $1$-form $a \in \mathcal{C}_{2+\nu}^{1,\gamma}(\mathbb{R}^8)$ 
such that 
\[\|a\|_{\mathcal{C}_{2+\nu}^{1,\gamma}(\mathbb{R}^8)} \leq C \, \|\psi\|_{\mathcal{C}_{3+
\nu}^\gamma(\mathbb{R}^8)}\] 
and 
\[\mathbb{L}_B a = \psi.\] 
Let $\zeta$ be a cut-off function on $E$ such that $\zeta(x) = 1$ for $|x| \leq \delta$, $\zeta
(x) = 0$ for $|x| \geq 2\delta$, and 
\[\sup \, \delta \, |\nabla \zeta| \leq C.\] 
Furthermore, let $\eta$ be a cut-off function on $E^\perp$ satisfying $\eta(y) = 1$ for $|y| 
\leq 2\delta^4$, $\eta(y) = 0$ for $|y| \geq 2\delta^2$, and 
\[\sup \, |y| \, |\nabla \eta| \leq C \, |\log \delta|^{-1}.\] 
Then we have the estimates 
\[\|\eta \, \zeta \, a\|_{\mathcal{C}_{2+\nu}^{1,\gamma}(\mathbb{R}^8)} \leq C \, \|\psi\|
_{\mathcal{C}_{3+\nu}^\gamma(\mathbb{R}^8)}\] 
and 
\begin{align*} 
&|\mathbb{L}_B(\zeta \, a) - \psi\|_{\mathcal{C}_{3+\nu}^\gamma(\{(x,y) \in \mathbb{R}^8: |y| 
\leq 2\delta^4\})} \\ 
&= \|\mathbb{L}_B (\zeta \, a) - \zeta \, \mathbb{L}_B a\|_{\mathcal{C}_{3+\nu}^\gamma(\{(x,y) 
\in \mathbb{R}^8: |y| \leq 2\delta^4\})} \\ 
&\leq C \, \delta \, \|a\|_{\mathcal{C}_{2+\nu}^\gamma(\mathbb{R}^8)} \\ 
&\leq C \, \delta \, \|\psi\|_{\mathcal{C}_{3+\nu}^\gamma(\mathbb{R}^8)} 
\end{align*} 
and 
\begin{align*} 
&|\mathbb{L}_B(\eta \, \zeta \, a) - \psi\|_{\mathcal{C}_{3+\nu}^\gamma(\mathbb{R}^8)} \\ 
&= \|\mathbb{L}_B (\eta \, \zeta \, a) - \eta \, \zeta \, \mathbb{L}_B a\|_{\mathcal{C}_{3+\nu}
^\gamma(\mathbb{R}^8)} \\ 
&\leq C \, |\log \delta|^{-1} \, \|a\|_{\mathcal{C}_{2+\nu}^\gamma(\mathbb{R}^8)} \\ 
&\leq C \, |\log \delta|^{-1} \, \|\psi\|_{\mathcal{C}_{3+\nu}^\gamma(\mathbb{R}^8)}. 
\end{align*} 
From this the assertion follows. \\

\section{Construction of the approximate solutions}

In this section, we outline the construction of a certain class of approximate solutions. To 
this end, we assume that the normal bundle $NS$ can be endowed with a $SU(2)$-structure 
$(J,\omega)$. Here, $J$ is a complex structure and $\omega$ is a complex volume form on $NS$. \\

Let $\nabla' = \nabla + \theta$ be a connection on the normal bundle $NS$ such that $\theta$ is 
a $1$-form with values in the Lie algebra $\Lambda_+^2 NS$ and $(J,\omega)$ is parallel with 
respect to the connection $\nabla'$. The $1$-form $\theta$ is uniquely determined by the 
covariant derivative of the pair $(J,\omega)$ with respect to the Levi-Civita connection 
$\nabla$. Since $(J,\omega)$ is parallel with respect to $\nabla'$, the connection induced by 
$\nabla'$ on the bundle $\Lambda_+^2 NS$ is flat. \\

The connection $\nabla'$ induces a splitting of the tangent space $TNS$ into horizontal and 
vertical subspaces. Let $\{e_i': 1 \leq i \leq 4\}$ be an orthonormal basis for the horizontal 
subspace with respect to $\nabla'$, and let $\{e_j^\perp: 1 \leq j \leq 4\}$ be a $SU(2)$ basis 
for the vertical subspace. \\

In the first step, we define a connection on the pull-back bundle $\pi^* NS$ of the normal 
bundle under the natural projection $\pi: NS \to S$. Since we may identify a neighborhood of $S$ 
in $M$ with a neighborhood of the zero section in $NS$, this gives a connection on a small 
neighborhood of $S$ in $M$. In the second step, we show that this connection can be extended to 
the whole of $M$ using suitable cut-off functions. \\

The glueing data consist of a set $(v,\lambda,J,\omega)$, where $v$ is a section of the normal 
bundle $NS$, $\lambda$ is a positive function on $S$, and $(J,\omega)$ is a $SU(2)$ structure on 
the normal bundle $NS$. Let $\{\mathfrak{i},\mathfrak{j},\mathfrak{k}\}$ be a basis for the Lie 
algebra $\mathfrak{su}(NS)$ such that 
\begin{align*} 
&\mathfrak{i}(e_1^\perp) = -e_2^\perp, \quad \mathfrak{i}(e_2^\perp) = e_1^\perp, \quad 
\mathfrak{i}(e_3^\perp) = e_4^\perp, \quad \mathfrak{i}(e_4^\perp) = -e_3^\perp, \\ 
&\mathfrak{j}(e_1^\perp) = -e_3^\perp, \quad \mathfrak{j}(e_2^\perp) = -e_4^\perp, \quad 
\mathfrak{j}(e_3^\perp) = e_1^\perp, \quad \mathfrak{j}(e_4^\perp) = e_2^\perp, \\ 
&\mathfrak{k}(e_1^\perp) = -e_4^\perp, \quad \mathfrak{k}(e_2^\perp) = e_3^\perp, \quad 
\mathfrak{k}(e_3^\perp) = -e_2^\perp, \quad \mathfrak{k}(e_4^\perp) = e_1^\perp. 
\end{align*} 
We consider a connection of the form $D_A = \nabla' + A$. The vertical components of $A$ are 
defined by 
\begin{align*} 
A(e_1^\perp) &= \frac{-(y - \varepsilon v)_2 \, \mathfrak{i} - (y - \varepsilon v)_3 \, 
\mathfrak{j} - (y - \varepsilon v)_4 \, \mathfrak{k}}{\varepsilon^2 \lambda^2 + |y - \varepsilon 
v|^2} \\ 
A(e_2^\perp) &= \frac{(y - \varepsilon v)_1 \, \mathfrak{i} - (y - \varepsilon v)_4 \, \mathfrak
{j} + (y - \varepsilon v)_3 \, \mathfrak{k}}{\varepsilon^2 \lambda^2 + |y - \varepsilon v|^2} \\ 
A(e_3^\perp) &= \frac{(y - \varepsilon v)_4 \, \mathfrak{i} + (y - \varepsilon v)_1 \, \mathfrak
{j} - (y - \varepsilon v)_2 \, \mathfrak{k}}{\varepsilon^2 \lambda^2 + |y - \varepsilon v|^2} \\ 
A(e_4^\perp) &= \frac{-(y - \varepsilon v)_3 \, \mathfrak{i} + (y - \varepsilon v)_2 \, 
\mathfrak{j} + (y - \varepsilon v)_1 \, \mathfrak{k}}{\varepsilon^2 \lambda^2 + |y - \varepsilon 
v|^2}. 
\end{align*} 
Since the basic instanton on $\mathbb{R}^4$ is $SU(2)$-equivariant, this definition is 
independent of the choice of $SU(2)$-frame $\{e_j^\perp: 1 \leq j \leq 4\}$. Furthermore, the 
horizontal components of $A$ are defined by 
\[A(e_i') = -\varepsilon \, \nabla_i' v_k \, A(e_k^\perp) - \lambda^{-1} \, \nabla_i \lambda \, 
(y - \varepsilon v)_k \, A(e_k^\perp)\] 
for $1 \leq i \leq 4$. \\

\begin{lemma}
The curvature of $A$ is given by 
\[F_A(e_i',e_j^\perp) = -\Big ( \varepsilon \, \nabla_i' v_k + \lambda^{-1} \, \nabla_i \lambda 
\, (y - \varepsilon v)_k \Big ) \, F_A(e_k^\perp,e_j^\perp)\] 
and 
\begin{align*} 
F_A(e_i',e_j') 
&= \Big ( \varepsilon \, \nabla_i' v_k + \lambda^{-1} \, \nabla_i \lambda \, (y - \varepsilon v)
_k \Big ) \\ 
&\cdot \Big ( \varepsilon \, \nabla_j' v_l + \lambda^{-1} \, \nabla_i \lambda \, (y - 
\varepsilon v)_l \Big ) \, F_A(e_k^\perp,e_l^\perp) \\ 
&+ C_{ij} + A \big ( C_{ij} \, (y - \varepsilon v) \big ), 
\end{align*} 
where $C_{ij} \in \Lambda_-^2 NS$ is the curvature of the connection $\nabla'$.
\end{lemma}

\vspace{2mm}

If $\{e_i: 1 \leq i \leq 4\}$ is an orthonormal basis for the horizontal subspace with respect 
to the Levi-Civita connection $\nabla$, then we obtain the following result: \\

\begin{lemma}
The curvature of $A$ satisfies 
\[F_A(e_i,e_j^\perp) = -\Big ( \varepsilon \, \nabla_i v_k + \lambda^{-1} \, \nabla_i \lambda \, 
(y - \varepsilon v)_k + \theta_{i,kl} \, (y - \varepsilon v)_l \Big ) \, F_A(e_k^\perp,e_j
^\perp)\] 
and 
\begin{align*} 
F_A(e_i,e_j) 
&= \Big ( \varepsilon \, \nabla_i v_k + \lambda^{-1} \, \nabla_i \lambda \, (y - \varepsilon v)
_k + \theta_{i,km} \, (y - \varepsilon v)_m \Big ) \\ 
&\cdot \Big ( \varepsilon \, \nabla_j v_l + \lambda^{-1} \, \nabla_j \lambda \, (y - \varepsilon 
v)_l + \theta_{j,ln} \, (y - \varepsilon v)_n \Big ) \, F_A(e_k^\perp,e_l^\perp) \\ 
&+ C_{ij} + A \big ( C_{ij} \, (y - \varepsilon v) \big ), 
\end{align*} 
where $C_{ij} \in \Lambda_-^2 NS$ is the curvature of $\nabla'$.
\end{lemma}

\vspace{2mm}

\begin{lemma}
Suppose that $\mu$ is constant and $r$ is a section of the vector bundle $\Lambda_+^2 NS$ such 
that $\nabla' r = 0$. Let 
\[u = (\varepsilon^2 \lambda^2  + |y - \varepsilon v|^2)^{-\frac{1}{2}} \, \big ( \mu \, (y - 
\varepsilon v)_k + r_{kl} \, (y - \varepsilon v)_l \big ) \, e_k^\perp.\] 
Then the covariant derivative of $u$ satisfies the estimate 
\[\|D_A u\|_{\mathcal{C}_3^\gamma(M)} \leq C \, \varepsilon^2.\] 
\end{lemma}

\vspace{2mm}

Hence, as we move away from the submanifold $S$, the connection $A$ approaches a flat 
connection. Therefore, we can extend $A$ trivially to $M$. \\

Our aim is to derive estimates for $F_A + *(\Omega \wedge F_A)$ in $\mathcal{C}_3^\gamma(M)$. To 
this end, we assume that the glueing data $(v,\lambda,J,\omega)$ satisfy the estimates 
\[\|v\|_{\mathcal{C}^{1,\gamma(M)}(S)} \leq K,\] 
\[\|\lambda\|_{\mathcal{C}^{1,\gamma(M)}(S)} \leq K, \qquad \inf \lambda \geq 1,\] 
\[\|(J,\omega)\|_{\mathcal{C}^{1,\gamma(M)}(S)} \leq K\] 
for some $K > 0$. All implicite constants will depend on $K$. \\

\begin{proposition}
If the set $(v,\lambda,J,\omega)$ is admissible, then we have the estimate 
\[\|F_A + *(\Omega \wedge F_A)\|_{\mathcal{C}_3^\gamma(M)} \leq C \, \varepsilon^2.\] 
\end{proposition}

\textbf{Proof.} 
Let $\Omega_0$ be a $4$-form which defines an almost $Spin(7)$-structure on $M$ such that 
$\Omega(x) = \Omega_0(x)$ for all $x \in S$ and $\nabla_X \Omega(x) = 0$ for all $x \in S$ and 
$X \in NS_x$. Then we have the estimate 
\begin{align*} 
&\|F_A + *(\Omega_0 \wedge F_A)\|_{\mathcal{C}_3^\gamma(M)} \\ 
&\leq \|F_A(e_1,e_2) + F_A(e_3,e_4) - F_A(e_1^\perp,e_2^\perp) - F_A(e_3^\perp,e_4^\perp)\|
_{\mathcal{C}_3^\gamma(M)} \\ 
&+ \|F_A(e_1,e_3) - F_A(e_2,e_4) - F_A(e_1^\perp,e_3^\perp) + F_A(e_2^\perp,e_4^\perp)\|
_{\mathcal{C}_3^\gamma(M)} \\ 
&+ \|F_A(e_1,e_4) + F_A(e_2,e_3) - F_A(e_1^\perp,e_4^\perp) - F_A(e_2^\perp,e_3^\perp)\|
_{\mathcal{C}_3^\gamma(M)} \\ 
&+ \|F_A(e_1,e_1^\perp) + F_A(e_2,e_2^\perp) + F_A(e_3,e_3^\perp) + F_A(e_4,e_4^\perp)\|
_{\mathcal{C}_3^\gamma(M)} \\ 
&+ \|F_A(e_1,e_2^\perp) - F_A(e_2,e_1^\perp) - F_A(e_3,e_4^\perp) + F_A(e_4,e_3^\perp)\|
_{\mathcal{C}_3^\gamma(M)} \\ 
&+ \|F_A(e_1,e_3^\perp) + F_A(e_2,e_4^\perp) - F_A(e_3,e_1^\perp) - F_A(e_4,e_2^\perp)\|
_{\mathcal{C}_3^\gamma(M)} \\ 
&+ \|F_A(e_1,e_4^\perp) - F_A(e_2,e_3^\perp) + F_A(e_3,e_2^\perp) - F_A(e_4,e_1^\perp)\|
_{\mathcal{C}_3^\gamma(M)}. 
\end{align*} 
Using the identities 
\begin{align*} 
&F_A(e_1^\perp,e_2^\perp) + F_A(e_3^\perp,e_4^\perp) = 0 \\ 
&F_A(e_1^\perp,e_3^\perp) + F_A(e_4^\perp,e_2^\perp) = 0 \\ 
&F_A(e_1^\perp,e_4^\perp) + F_A(e_2^\perp,e_3^\perp) = 0, 
\end{align*} 
we obtain 
\[\|F_A + *(\Omega_0 \wedge F_A)\|_{\mathcal{C}_3^\gamma(M)} \leq C \, \varepsilon^2.\] 
Since $\Omega = \Omega_0 + O(|y|)$, we conclude that 
\[\|F_A + *(\Omega \wedge F_A)\|_{\mathcal{C}_3^\gamma(M)} \leq \|F_A + *(\Omega_0 \wedge F_A)\|
_{\mathcal{C}_3^\gamma(M)} + \|F_A\|_{\mathcal{C}_4^\gamma(M)} \leq C \, \varepsilon^2.\] 

\vspace{2mm}

\section{Estimates for the linearized operator in weighted H\"older spaces}

Our aim in this section is to analyze the mapping properties of the linearized operator 
$\mathbb{L}_A: \Omega^1(M) \to \Omega_+^2(M)$. \\

As in Section 2, we define two vector bundles $V \subset TS \otimes NS$ and $W \subset TS 
\otimes NS \otimes NS$ over the submanifold $S$. Both vector bundles have rank $4$. The 
following elements form a basis for $V$: 
\[e_1 \otimes e_1^\perp + e_2 \otimes e_2^\perp + e_3 \otimes e_3^\perp + e_4 \otimes e_4^\perp,
\] 
\[e_1 \otimes e_2^\perp - e_2 \otimes e_1^\perp - e_3 \otimes e_4^\perp + e_4 \otimes e_3^\perp,
\] 
\[e_1 \otimes e_3^\perp + e_2 \otimes e_4^\perp - e_3 \otimes e_1^\perp - e_4 \otimes e_2^\perp,
\] 
\[e_1 \otimes e_4^\perp - e_2 \otimes e_3^\perp + e_3 \otimes e_2^\perp - e_4 \otimes e_1^\perp.
\] 
Similarly, the following elements form a basis for $W$: 
\begin{align*} 
&e_1 \otimes (e_1^\perp \otimes e_1^\perp + e_2^\perp \otimes e_2^\perp + e_3^\perp \otimes e
_3^\perp + e_4^\perp \otimes e_4^\perp) \\ 
&+ e_2 \otimes (e_2^\perp \otimes e_1^\perp - e_1^\perp \otimes e_2^\perp + e_4^\perp \otimes e
_3^\perp - e_3^\perp \otimes e_4^\perp) \\ 
&+ e_3 \otimes (e_3^\perp \otimes e_1^\perp - e_4^\perp \otimes e_2^\perp - e_1^\perp \otimes e
_3^\perp + e_2^\perp \otimes e_4^\perp) \\ 
&+ e_4 \otimes (e_4^\perp \otimes e_1^\perp + e_3^\perp \otimes e_2^\perp - e_2^\perp \otimes e
_3^\perp - e_1^\perp \otimes e_4^\perp), 
\end{align*} 
\begin{align*} 
&e_1 \otimes (e_2^\perp \otimes e_1^\perp - e_1^\perp \otimes e_2^\perp + e_4^\perp \otimes e
_3^\perp - e_3^\perp \otimes e_4^\perp) \\ 
&- e_2 \otimes (e_1^\perp \otimes e_1^\perp + e_2^\perp \otimes e_2^\perp + e_3^\perp \otimes e
_3^\perp + e_4^\perp \otimes e_4^\perp) \\ 
&- e_3 \otimes (e_4^\perp \otimes e_1^\perp + e_3^\perp \otimes e_2^\perp - e_2^\perp \otimes e
_3^\perp - e_1^\perp \otimes e_4^\perp) \\  
&+ e_4 \otimes (e_3^\perp \otimes e_1^\perp - e_4^\perp \otimes e_2^\perp - e_1^\perp \otimes e
_3^\perp + e_2^\perp \otimes e_4^\perp), 
\end{align*}
\begin{align*} 
&e_1 \otimes (e_3^\perp \otimes e_1^\perp - e_4^\perp \otimes e_2^\perp - e_1^\perp \otimes e
_3^\perp + e_2^\perp \otimes e_4^\perp) \\ 
&+ e_2 \otimes (e_4^\perp \otimes e_1^\perp + e_3^\perp \otimes e_2^\perp - e_2^\perp \otimes e
_3^\perp - e_1^\perp \otimes e_4^\perp) \\ 
&- e_3 \otimes (e_1^\perp \otimes e_1^\perp + e_2^\perp \otimes e_2^\perp + e_3^\perp \otimes e
_3^\perp + e_4^\perp \otimes e_4^\perp) \\ 
&- e_4 \otimes (e_2^\perp \otimes e_1^\perp - e_1^\perp \otimes e_2^\perp + e_4^\perp \otimes e
_3^\perp - e_3^\perp \otimes e_4^\perp), 
\end{align*} 
\begin{align*}
&e_1 \otimes (e_4^\perp \otimes e_1^\perp + e_3^\perp \otimes e_2^\perp - e_2^\perp \otimes e
_3^\perp - e_1^\perp \otimes e_4^\perp) \\ 
&- e_2 \otimes (e_3^\perp \otimes e_1^\perp - e_4^\perp \otimes e_2^\perp - e_1^\perp \otimes e
_3^\perp + e_2^\perp \otimes e_4^\perp) \\ 
&+ e_3 \otimes (e_2^\perp \otimes e_1^\perp - e_1^\perp \otimes e_2^\perp + e_4^\perp \otimes e
_3^\perp - e_3^\perp \otimes e_4^\perp) \\ 
&- e_4 \otimes (e_1^\perp \otimes e_1^\perp + e_2^\perp \otimes e_2^\perp + e_3^\perp \otimes e
_3^\perp + e_4^\perp \otimes e_4^\perp). 
\end{align*}

\begin{proposition}
Suppose that $\psi \in \mathcal{C}_{3+\nu}^\gamma(M)$ is supported in the set $\{p \in M: \dist
(p,S) \leq 2\delta^4\}$ and satisfies 
\[\int_{NS_x} \sum_{i,j=1}^4 \big ( \varepsilon \, s_{ik} + t_{ikl} \, (y - \varepsilon v)_l 
\big ) \, \langle \psi(e_i,e_j^\perp),F_A(e_k^\perp,e_j^\perp) \rangle = 0\] 
for all $x \in S$, $s \in V_x$, and $t \in W_x$. Then there exists a $1$-form $a \in \mathcal{C}
_{2+\nu}^{1,\gamma}(M)$ which is supported in the region $\{p \in M: \dist(p,S) \leq 2\delta^2
\}$ such that 
\[\|a\|_{\mathcal{C}_{2+\nu}^{1,\gamma}(M)} \leq C \, \|\psi\|_{\mathcal{C}_{3+\nu}^\gamma(M)}\] 
and 
\[\|\mathbb{L}_A a - \psi\|_{\mathcal{C}_{3+\nu}^\gamma(\{p \in M: \dist(p,S) \leq 2\delta^4\})} 
\leq C \, \delta \, \|\psi\|_{\mathcal{C}_{3+\nu}^\gamma(M)},\] 
and 
\[\|\mathbb{L}_A a - \psi\|_{\mathcal{C}_{3+\nu}^\gamma(M)} \leq C \, |\log \delta|^{-1} \, 
\|\psi\|_{\mathcal{C}_{3+\nu}^\gamma(M)}.\] 
\end{proposition}

\textbf{Proof.} 
Let $\{\zeta^{(j)}: 1 \leq j \leq j_0\}$ be a partition of unity on $S$ such that each function 
$\zeta^{(j)}$ is supported in a ball $B_\delta(p_j)$, and 
\[|\{1 \leq j \leq j_0: x \in B_{4\delta}(p_j)\}| \leq C\] 
for all $x \in S$ and some uniform constant $C$. For each $1 \leq j \leq j_0$, there exists a 
$1$-form $a^{(j)} \in \mathcal{C}_{2+\nu}^{1,\gamma}(M)$ which is supported in the region 
$\{(x,y) \in NS: x \in B_{2\delta}(p_j), \, |y| \leq 2\delta^2\}$ such that 
\[\|a^{(j)}\|_{\mathcal{C}_{2+\nu}^{1,\gamma}(M)} \leq C \, \|\zeta^{(j)} \, \psi\|_{\mathcal{C}
_{3+\nu}^\gamma(M)}\] 
and 
\[\|\mathbb{L}_A a^{(j)} - \zeta^{(j)} \, \psi\|_{\mathcal{C}_{3+\nu}^\gamma(\{p \in M: \dist(p,
S) \leq 2\delta^4\})} \leq C \, \delta \, \|\zeta^{(j)} \, \psi\|_{\mathcal{C}_{3+\nu}^\gamma
(M)},\] 
and 
\[\|\mathbb{L}_A a^{(j)} - \zeta^{(j)} \, \psi\|_{\mathcal{C}_{3+\nu}^\gamma(M)} \leq C \, |\log 
\delta|^{-1} \, \|\zeta^{(j)} \, \psi\|_{\mathcal{C}_{3+\nu}^\gamma(M)}.\] 
We now define 
\[a = \sum_{j=1}^{j_0} a^{(j)}.\] 
Then we have the estimates 
\begin{align*} 
\|a\|_{\mathcal{C}_{2+\nu}^{1,\gamma}(M)} 
&\leq C \, \sup_{1 \leq j \leq j_0} \|a^{(j)}\|_{\mathcal{C}_{2+\nu}^{1,\gamma}(M)} \\ 
&\leq C \, \sup_{1 \leq j \leq j_0} \|\zeta^{(j)} \, \psi\|_{\mathcal{C}_{3+\nu}^\gamma(M)} \\ 
&\leq C \, \|\psi\|_{\mathcal{C}_{3+\nu}^\gamma(M)}, 
\end{align*} 
\begin{align*} 
\|\mathbb{L}_A a - \psi\|_{\mathcal{C}_{3+\nu}^\gamma(\{p \in M: \dist(p,S) \leq 2\delta^4\})} 
&\leq C \, \sup_{1 \leq j \leq j_0} \|\mathbb{L}_A a^{(j)} - \zeta^{(j)} \, \psi\|_{\mathcal{C}
_{3+\nu}^\gamma(\{p \in M: \dist(p,S) \leq 2\delta^4\})} \\ 
&\leq C \, \delta \, \sup_{1 \leq j \leq j_0} \|\zeta^{(j)} \, \psi\|_{\mathcal{C}_{3+\nu}
^\gamma(M)} \\ 
&\leq C \, \delta \, \|\psi\|_{\mathcal{C}_{3+\nu}^\gamma(M)}, 
\end{align*} 
\begin{align*} 
\|\mathbb{L}_A a - \psi\|_{\mathcal{C}_{3+\nu}^\gamma(M)} 
&\leq C \, \sup_{1 \leq j \leq j_0} \|\mathbb{L}_A a^{(j)} - \zeta^{(j)} \, \psi\|_{\mathcal{C}
_{3+\nu}^\gamma(M)} \\ 
&\leq C \, |\log \delta|^{-1} \, \sup_{1 \leq j \leq j_0} \|\zeta^{(j)} \, \psi\|_{\mathcal{C}
_{3+\nu}^\gamma(M)} \\ 
&\leq C \, |\log \delta|^{-1} \, \|\psi\|_{\mathcal{C}_{3+\nu}^\gamma(M)}. 
\end{align*} 
This proves the assertion. \\

\begin{proposition}
For every $\Omega$-self-dual $2$-form $\psi \in \mathcal{C}_{3+\nu}^\gamma(M)$, there exists a 
$1$-form $a \in \mathcal{C}_{2+\nu}^{1,\gamma}(M)$ such that 
\[\|a\|_{\mathcal{C}_{2+\nu}^{1,\gamma}(M)} \leq C \, \|\psi\|_{\mathcal{C}_{3+\nu}^\gamma(M)}\] 
and 
\[P_+ da = \psi.\]
\end{proposition}

\textbf{Proof.} 
We consider the elliptic operator $P_+ dd^*: \Omega_+^2(M) \to \Omega_+^2(M)$. Its kernel is 
given by 
\[\ker (P_+ dd^*: \Omega_+^2(M) \to \Omega_+^2(M)) = H_+^2(M).\] 
Since the cohomology group $H_+^2(M)$ vanishes, the operator $P_+ dd^*: \Omega_+^2(M) \to \Omega
_+^2(M)$ is invertible. Consequently, there exists a $\Omega$-self-dual $2$-form $\varphi$ such 
that 
\[P_+ dd^* \varphi = \psi.\] 
We claim that 
\[\|\varphi\|_{\mathcal{C}_{1+\nu}^{2,\gamma}(M)} \leq C \, \|P_+ dd^* \varphi\|_{\mathcal{C}_{3
+\nu}^\gamma(M)}.\] 
By Schauder estimates, it suffices to show that 
\[\sup \, (\varepsilon + \dist(p,S))^{1+\nu} \, |\varphi| \leq C \, \sup \, (\varepsilon + \dist
(p,S))^{3+\nu} \, |P_+ dd^* \varphi|.\] 
If this estimate fails, then there exists a sequence of positive real numbers $\varepsilon_j$ 
and a sequence of $\Omega$-self-dual $2$-forms $\varphi^{(j)} \in \mathcal{C}_{1+\nu}^{2,\gamma}
(M)$ such that 
\[\sup \, (\varepsilon_j + \dist(p,S))^{1+\nu} \, |\varphi^{(j)}| = 1\] 
and 
\[\sup \, (\varepsilon_j + \dist(p,S))^{3+\nu} \, |P_+ dd^* \varphi^{(j)}| \to 0.\] 
Then there exists a sequence of points $p_j \in M$ such 
that 
\[\sup \, (\varepsilon_j + \dist(p_j,S))^{1+\nu} \, |\varphi^{(j)}(p_j)| \geq \frac{1}{2}.\] 
There are two possibilities: \\ 

(i) Suppose that $\dist(p_j,S)$ is bounded from below. After passing to a subsequence, we may 
assume that the sequence $\varphi^{(j)}$ converges to a $\Omega$-self-dual $2$-form $\varphi \in 
\Omega_+^2(M)$ such that 
\[\sup \, \dist(p,S)^{1+\nu} \, |\varphi| \leq 1\] 
and 
\[P_+ dd^* \varphi = 0.\] 
From this it follows that $\varphi$ is smooth. Since the operator $P_+ dd^*: \Omega_+^2(M) \to 
\Omega_+^2(M)$ has trivial kernel, it follows that $\varphi = 0$. This is a contradiction. \\

(ii) We now assume that $\dist(p_j,S) \to 0$. After rescaling and taking the limit, we obtain a 
$\Omega$-self-dual $2$-form $\tilde{\varphi} \in \Omega_+^2(\mathbb{R}^8)$ such that 
\[\sup \, |y|^{1+\nu} \, |\tilde{\varphi}| \leq 1\] 
and 
\[P_+ dd^* \tilde{\varphi} = 0.\] 
Thus, we conclude that $\tilde{\varphi} = 0$. This is a contradiction. \\

This implies 
\[\|\varphi\|_{\mathcal{C}_{1+\nu}^{2,\gamma}(M)} \leq C \, \|\psi\|_{\mathcal{C}_{3+\nu}^\gamma
(M)}.\] 
Letting $a = d^* \varphi$, the assertion follows. \\

\begin{proposition}
Suppose that $\psi \in \mathcal{C}_{3+\nu}^\gamma(M)$ is supported in the region $\{p \in M: 
\dist(p,S) \geq \delta^4\}$. Then there exists a $1$-form $a \in \mathcal{C}_{2+\nu}^{1,\gamma}
(M)$ which is supported in the region $\{p \in M: \dist(p,S) \geq \delta^8\}$ such that 
\[\|a\|_{\mathcal{C}_{2+\nu}^{1,\gamma}(M)} \leq C \, \|\psi\|_{\mathcal{C}_{3+\nu}^\gamma(M)}\] 
and 
\[\|\mathbb{L}_A a - \psi\|_{\mathcal{C}_{3+\nu}^\gamma(M)} \leq C \, \Big ( |\log \delta|^{-1} 
+ \delta^{-16} \, \varepsilon^2 \Big ) \, \|\psi\|_{\mathcal{C}_{3+\nu}^\gamma(M)}.\]
\end{proposition}

\textbf{Proof.} 
By Proposition 4.2, exists a $1$-form $a$ such that 
\[\|a\|_{\mathcal{C}_{2+\nu}^{1,\gamma}(M)} \leq C \, \|\psi\|_{\mathcal{C}_{3+\nu}^\gamma(M)}\] 
and 
\[2 \, P_+ da = \psi.\] 
Let $\eta$ be a cut-off function such that $\eta(p) = 0$ for $\dist(p,S) \leq \delta^8$, $\eta
(p) = 1$ for $\dist(p,S) \geq \delta^4$ and 
\[\sup \, \dist(p,S) \, |\nabla \eta| \leq C \, |\log \delta|^{-1}.\] 
Then the $1$-form $\eta \, a$ is supported in the region $\{p \in M: \dist(p,S) \geq \delta^8\}$ 
and satisfies 
\begin{align*} 
&\|\mathbb{L}_A (\eta \, a) - \psi\|_{\mathcal{C}_{3+\nu}^\gamma(M)} \\ 
&= 2 \, \|P_+ D_A(\eta \, a)  - \eta \, P_+ da\|_{\mathcal{C}_{3+\nu}^\gamma(M)} \\ 
&\leq 2 \, \|P_+ D_A(\eta \, a) - P_+ d(\eta \, a)\|_{\mathcal{C}_{3+\nu}^\gamma(M)} + 2 \, \|P
_+ d(\eta \, a) - \eta \, P_+ da\|_{\mathcal{C}_{3+\nu}^\gamma(M)} \\ 
&\leq C \, \delta^{-16} \, \varepsilon^2 \, \|a\|_{\mathcal{C}_{2+\nu}^\gamma(M)} + C \, |\log 
\delta|^{-1} \, \|a\|_{\mathcal{C}_{2+\nu}^\gamma(M)} \\ 
&\leq C \, \delta^{-16} \, \varepsilon^2 \, \|\psi\|_{\mathcal{C}_{3+\nu}^\gamma(M)} + C \, 
|\log \delta|^{-1} \, \|\psi\|_{\mathcal{C}_{3+\nu}^\gamma(M)}. 
\end{align*} 
This proves the assertion. \\

In the following, we will choose $\delta = \varepsilon^{\frac{1}{16}}$. Let $\kappa$ be a 
cut-off function such that $\kappa(p) = 1$ for $\dist(p,S) \leq \varepsilon^{\frac{1}{4}}$ and 
$\kappa(p) = 0$ for $\dist(p,S) \geq 2 \, \varepsilon^{\frac{1}{4}}$. \\

Let $\mathcal{G}_\nu^{k,\gamma}(M)$ be the set of all $\psi \in \Omega_+^2(M)$ such that $\psi 
\in \mathcal{C}_\nu^{k,\gamma}(M)$ and 
\[\int_{NS_x} \kappa \, \sum_{i,j=1}^4 \big ( \varepsilon \, s_{ik} + t_{ikl} \, (y - 
\varepsilon v)_l \big ) \, \langle \psi(e_i,e_j^\perp),F_A(e_k^\perp,e_j^\perp) \rangle = 0\] 
for all $x \in S$, $s \in V_x$, and $t \in W_x$. \\

We denote by $\Id - \mathbb{P}$ the fibrewise projection from $\mathcal{C}_\nu^\gamma(M)$ to the 
subspace $\mathcal{G}_\nu^\gamma(M)$. Hence, if $\psi$ is an $\Omega$-self-dual $2$-form, then 
the projection $\mathbb{P} \psi$ is of the form 
\[\mathbb{P} \psi(e_i,e_j^\perp) = \kappa \, \big ( \varepsilon \, s_{ik} + t_{ikl} \, (y - 
\varepsilon v)_l \big ) \, F_A(e_k^\perp,e_j^\perp)\] 
for suitable $s \in V$ and $t \in W$. Let $\Pi$ be the linear operator which assigns to every 
$\Omega$-self-dual $2$-form $\psi$ the pair 
\[\Pi \psi = (s,t) \in V \oplus W.\] 
We shall need the following estimate for the operator norm of the projection operator $\mathbb
{P}$. \\

\begin{proposition}
For every $\Omega$-self-dual $2$-form $\psi \in \mathcal{C}_{3+\nu}^\gamma(M)$, we have the 
estimates 
\[\|\Pi \psi\|_{\mathcal{C}^\gamma(S)} \leq C \, \varepsilon^{-2-\nu-\gamma} \, \|\psi\|
_{\mathcal{C}_{3+\nu}^\gamma(M)}\] 
and 
\[\|\mathbb{P} \psi\|_{\mathcal{C}_{3+\nu}^\gamma(M)} \leq C \, \varepsilon^{-\nu-\gamma} \, 
\|\psi\|_{\mathcal{C}_{3+\nu}^\gamma(M)}.\] 
\end{proposition}

\textbf{Proof.} 
This follows from \cite{Br}, Proposition 5.4. \\

\begin{proposition}
For every $\psi \in \mathcal{G}_{3+\nu}^\gamma(M)$ there exists a $1$-form $a \in \mathcal{C}
_{2+\nu}^{1,\gamma}(M)$ such that 
\[\|a\|_{\mathcal{C}_{2+\nu}^{1,\gamma}(M)} \leq C \, \|\psi\|_{\mathcal{C}_{3+
\nu}^\gamma(M)}\] 
and 
\[\|\mathbb{L}_A a - \psi\|_{\mathcal{C}_{3+\nu}^\gamma(\{p \in M: \dist(p,S) \leq \varepsilon
^{\frac{1}{2}}\})} \leq C \, \varepsilon^{\frac{1}{16}} \, \|\psi\|
_{\mathcal{C}_{3+\nu}^\gamma(M)},\] 
and 
\[\|\mathbb{L}_A a - \psi\|_{\mathcal{C}_{3+\nu}^\gamma(M)} \leq C \, |\log \varepsilon|^{-1} \, 
\|\psi\|_{\mathcal{C}_{3+\nu}^\gamma(M)}.\] 
\end{proposition}

\textbf{Proof.} 
Apply Proposition 4.1 to $\kappa \, \psi$ and Proposition 4.3 to $(1 - \kappa) \, \psi$. \\

\begin{proposition}
For every $\psi \in \mathcal{G}_{3+\nu}^\gamma(M)$ there exists a $1$-form $a \in \mathcal{C}
_{2+\nu}^{1,\gamma}(M)$ such that 
\[\|a\|_{\mathcal{C}_{2+\nu}^{1,\gamma}(M)} \leq C \, \|\psi\|_{\mathcal{C}_{3+\nu}^\gamma(M)}\] 
and 
\[(\Id - \mathbb{P}) \, \mathbb{L}_A a = \psi.\] 
Furthermore, $a$ satisfies the estimate 
\[\|\Pi \, \mathbb{L}_A a\|_{\mathcal{C}^\gamma(S)} \leq C \, \varepsilon^{-2+\frac{1}{32}} \, 
\|\psi\|_{\mathcal{C}_{3+\nu}^\gamma(M)}.\] 
\end{proposition}

\textbf{Proof.} 
By Proposition 4.5, there exists an operator $\mathbb{S}: \mathcal{G}_{3+\nu}^\gamma(M) \to 
\mathcal{C}_{2+\nu}^{1,\gamma}(M)$ such that 
\[\|\mathbb{S} \psi\|_{\mathcal{C}_{2+\nu}^{1,\gamma}(M)} \leq C \, \|\psi\|_{\mathcal{C}_{3+
\nu}^\gamma(M)}\] 
and 
\[\|\mathbb{L}_A \, \mathbb{S} \psi - \psi\|_{\mathcal{C}_{3+\nu}^\gamma(\{p \in M: \dist(p,S) 
\leq \varepsilon^{\frac{1}{2}}\})} \leq C \, \varepsilon^{\frac{1}{16}} \, \|\psi\|_{\mathcal{C}
_{3+\nu}^\gamma(M)},\] 
and 
\[\|\mathbb{L}_A \, \mathbb{S} \psi - \psi\|_{\mathcal{C}_{3+\nu}^\gamma(M)} \leq C \, |\log 
\varepsilon|^{-1} \, \|\psi\|_{\mathcal{C}_{3+\nu}^\gamma(M)}.\] 
This implies 
\[\|\Pi \, \mathbb{L}_A \, \mathbb{S} \psi\|_{\mathcal{C}^\gamma(S)} = \|\Pi (\mathbb{L}_A \, 
\mathbb{S} \psi - \psi)\|_{\mathcal{C}^\gamma(S)} \leq C \, \varepsilon^{-2+\frac{1}{16}-\nu-
\gamma} \, \|\psi\|_{\mathcal{C}_{3+\nu}^\gamma(M)}.\] 
From this it follows that 
\[\|(\Id - \mathbb{P}) \, \mathbb{L}_A \, \mathbb{S} \psi - \psi\|_{\mathcal{C}_{3+\nu}^\gamma
(M)} \leq C \, |\log \varepsilon|^{-1} \, \|\psi\|_{\mathcal{C}_{3+\nu}^\gamma(M)}.\] 
Therefore, the operator $(\Id - \mathbb{P}) \, \mathbb{L}_A \, \mathbb{S}: \mathcal{G}_{3+\nu}
^\gamma(M) \to \mathcal{G}_{3+\nu}^\gamma(M)$ is invertible. Hence, if we define 
\[a = \mathbb{S} \, \big [ (\Id - \mathbb{P}) \, \mathbb{L}_A \, \mathbb{S} \big ]^{-1} \, \psi,
\] 
then $a$ satisfies 
\[\|a\|_{\mathcal{C}_{2+\nu}^{1,\gamma}(M)} \leq C \, \|\psi\|_{\mathcal{C}_{3+\nu}^\gamma(M)}\] 
and 
\[(\Id - \mathbb{P}) \, \mathbb{L}_A a = \psi.\] 
This proves the assertion. \\

\section{The nonlinear problem}

\begin{proposition}
For every approximate solution $A$, there exists a nearby connection $\tilde{A} = A + a$ such 
that 
\[\|a\|_{C_{2+\nu}^{1,\gamma}(M)} \leq C \, \varepsilon^{2-\nu-\gamma}\] 
and 
\[(\Id - \mathbb{P}) \, (F_{\tilde{A}} + *(\Omega \wedge F_{\tilde{A}})) = 0.\] 
Furthermore, $a$ satisfies the estimate 
\[\|\Pi \, \mathbb{L}_A a\|_{\mathcal{C}^\gamma(S)} \leq C \, \varepsilon^{\frac{1}{32}}.\] 
\end{proposition}

\textbf{Proof.} 
The connection $\tilde{A} = A + a$ satisfies 
\[F_{\tilde{A}} + *(\Omega \wedge F_{\tilde{A}}) = F_A + *(\Omega \wedge F_A) + D_A a + *(\Omega 
\wedge D_A a) + [a,a] + *(\Omega \wedge [a,a]).\] 
This implies 
\[F_{\tilde{A}} + *(\Omega \wedge F_{\tilde{A}}) = F_A + *(\Omega \wedge F_A) + 2 \, \mathbb{L}
_A a + [a,a] + *(\Omega \wedge [a,a]).\] 
According to Proposition 4.6, there exists an operator $\mathbb{G}: \mathcal{G}_{3+\nu}^\gamma
(M) \to \mathcal{C}_{2+\nu}^{1,\gamma}(M)$ such that 
\[\|\mathbb{G} \psi\|_{\mathcal{C}_{2+\nu}^{1,\gamma}(M)} \leq C \, \|\psi\|_{\mathcal{C}_{3+
\nu}^\gamma(M)}\] 
and 
\[(\Id - \mathbb{P}) \, \mathbb{L}_A \, \mathbb{G} = \Id.\] 
We now define a mapping $\Phi: \mathcal{C}_{2+\nu}^{1,\gamma}(M) \to \mathcal{C}_{2+\nu}^{1,
\gamma}(M)$ by 
\[\Phi(a) = -\frac{1}{2} \, \mathbb{G} \, (\Id - \mathbb{P}) \, \big ( F_A + *(\Omega \wedge 
F_A) \big ) - \frac{1}{2} \, \mathbb{G} \, (\Id - \mathbb{P}) \, \big ( [a,a] + *(\Omega \wedge 
[a,a]) \big ).\] 
Then we have the estimate 
\begin{align*} 
\|\Phi(a)\|_{\mathcal{C}_{2+\nu}^{1,\gamma}(M)} 
&\leq C \, \big \| (\Id - \mathbb{P}) \, \big ( F_A + *(\Omega \wedge F_A) \big ) \big \|
_{\mathcal{C}_{3+\nu}^\gamma(M)} \\ 
&+ C \, \big \| (\Id - \mathbb{P}) \, \big ( [a,a] + *(\Omega \wedge [a,a]) \big ) \big \|
_{\mathcal{C}_{3+\nu}^\gamma(M)} \\ 
&\leq C \, \varepsilon^{-\nu-\gamma} \, \|F_A + *(\Omega \wedge F_A)\|_{\mathcal{C}_{3+\nu}
^\gamma(M)} \\ 
&+ C \, \varepsilon^{-\nu-\gamma} \, \|[a,a]\|_{\mathcal{C}_{3+\nu}
^\gamma(M)} \\ 
&\leq C \, \varepsilon^{-\nu-\gamma} \, \|F_A + *(\Omega \wedge F_A)\|_{\mathcal{C}_{3+\nu}
^\gamma(M)} \\ 
&+ C \, \varepsilon^{-1-2\nu-\gamma} \, \|a\|_{\mathcal{C}_{2+\nu}^{1,\gamma}(M)}^2 \\ 
&\leq C \, \varepsilon^{2-\nu-\gamma} 
\end{align*} 
for all $a \in \mathcal{C}_{2+\nu}^{1,\gamma}(M)$ satisfying 
\[\|a\|_{\mathcal{C}_{2+\nu}^{1,\gamma}(M)} \leq \varepsilon^{\frac{7}{4}}.\] 
Moreover, we have 
\begin{align*} \|\Phi(a) - \Phi(a')\|_{\mathcal{C}_{1+\nu}^{2,\gamma}(M)} 
&\leq C \, \varepsilon^{-\nu-\gamma} \, \|[a,a] - [a',a']\|_{\mathcal{C}_{3+\nu}^\gamma(M)} \\ 
&\leq C \, \varepsilon^{\frac{3}{4}-2\nu-\gamma} \, \|a - a'\|_{\mathcal{C}_{1+\nu}^{2,\gamma}
(M)} 
\end{align*} 
for all $a,a' \in \mathcal{C}_{2+\nu}^{1,\gamma}(M)$ satisfying 
\[\|a\|_{\mathcal{C}_{2+\nu}^{1,\gamma}(M)}, \, 
\|a'\|_{\mathcal{C}_{2+\nu}^{1,\gamma}(M)} \leq \varepsilon^{\frac{7}{4}}.\] 
Hence, it follows from the contraction mapping principle that there exists a $1$-form $a \in 
\mathcal{C}_{2+\nu}^{1,\gamma}(M)$ such that 
\[\|a\|_{\mathcal{C}_{2+\nu}^{1,\gamma}(M)} \leq C \, \varepsilon^{2-\nu-\gamma}\] 
and 
\[\Phi(a) = a.\] 
From this it follows that 
\[\mathbb{G} \, (\Id - \mathbb{P}) \, \big ( F_A + *(\Omega \wedge F_A) \big ) + 2a + \mathbb{G} 
\, (\Id - \mathbb{P}) \, \big ( [a,a] + *(\Omega \wedge [a,a] \big ) = 0,\] 
hence 
\[(\Id - \mathbb{P}) \, \big ( F_A + *(\Omega \wedge F_A) \big ) + 2 \, (\Id - \mathbb{P}) \, 
\mathbb{L}_A a + (\Id - \mathbb{P}) \, \big ( [a,a] + *(\Omega \wedge [a,a] \big ) = 0.\] 
Thus, we conclude that 
\[(\Id - \mathbb{P}) \, \big ( F_{\tilde{A}} + *(\Omega \wedge F_{\tilde{A}}) \big ) = 0.\] 
This proves the assertion. \\

\begin{corollary}
If $\tilde{A}$ satisfies 
\[\mathbb{P} \, (F_{\tilde{A}} + *(\Omega \wedge F_{\tilde{A}})) = 0,\] 
then $\tilde{A}$ is an $\Omega$-anti-self-dual instanton, i.e. 
\[F_{\tilde{A}} + *(\Omega \wedge F_{\tilde{A}}) = 0.\]
\end{corollary}

\vspace{2mm}

\section{The balancing condition}

By Corollary 5.2, the problem is reduced to finding a set of glueing data $(v,\lambda,J,\omega)$ 
such that 
\[\mathbb{P} \, (F_{\tilde{A}} + *(\Omega \wedge F_{\tilde{A}})) = 0.\] 
Our aim in this section is to derive a formula for the error term 
\[\mathbb{P} \, (F_{\tilde{A}} + *(\Omega \wedge F_{\tilde{A}})).\] 

\vspace{2mm}

\begin{proposition} 
The curvature of $A$ satisfies 
\begin{align*} 
\Pi(F_A + *(\Omega_0 \wedge F_A)) = 4 \, \bigg ( 
&\proj_V \Big ( \sum_{i,j=1}^4 \nabla_i v_k \, e_i \otimes e_k^\perp \Big ), \\ 
&\proj_W \Big ( \sum_{i,k,l=1}^4 (\lambda^{-1} \, \nabla_i \lambda \, \delta_{kl} + \theta
_{i,kl}) \, e_i \otimes e_k^\perp \otimes e_l^\perp \bigg ). 
\end{align*}
\end{proposition}

\textbf{Proof.} 
This is a consequence of the identity 
\[F_A(e_i,e_j^\perp) = -\Big ( \varepsilon \, \nabla_i v_k + \lambda^{-1} \, \nabla_i \lambda \, 
(y - \varepsilon v)_k + \theta_{i,kl} \, (y - \varepsilon v)_l \Big ) \, F_A(e_k^\perp,e_j
^\perp).\]

\vspace{2mm}

The covariant derivative of $\Omega$ can be described by a $1$-form $\alpha$ with values in 
$\Lambda_+^2 TM$. For every vector field $X \in TM$, we write 
\[\nabla_X \Omega = \sum_{k=1}^8 i_{e_k} \alpha(X) \wedge i_{e_k} \Omega,\] 
where $\alpha(X) \in \Lambda_+^2 TM$. From this it follows that 
\[\Omega = \Omega_0 + \sum_{k=1}^8 i_{e_k} \alpha(y) \wedge i_{e_k} \Omega_0 + O(|y|^2),\] 
where $\alpha(y) \in \Lambda_+^2 TM$. \\

\begin{proposition}
The curvature of $A$ satisfies 
\begin{align*} 
\bigg \| &\Pi(F_A + *(\Omega \wedge F_A)) \\ 
&- 4 \, \bigg ( \proj_V \Big ( \sum_{i,j=1}^4 (\nabla_i v_k + \alpha_{ik,l} \, v_l) \, e_i 
\otimes e_k^\perp \Big ), \\ 
&\hspace{10.1mm} \proj_W \Big ( \sum_{i,k,l=1}^4 (\lambda^{-1} \, \nabla_i \lambda \, \delta
_{kl} + \theta_{i,kl} + \alpha_{ik,l}) \, e_i \otimes e_k^\perp \otimes e_l^\perp \bigg ) \bigg 
\|_{\mathcal{C}^\gamma(S)} \leq C \, \varepsilon. 
\end{align*}
\end{proposition}

\textbf{Proof.} 
Using the identity 
\[\Omega - \Omega_0 - \sum_{k=1}^8 i_{e_k} \alpha(y) \wedge i_{e_k} \Omega_0 = O(|y|^2),\] 
we obtain 
\[\bigg \| \Omega \wedge F_A - \Omega_0 \wedge F_A - \sum_{k=1}^8 i_{e_k} \alpha(y) \wedge i_{e
_k} \Omega_0 \wedge F_A \bigg \|_{\mathcal{C}_2^\gamma(M)} \leq C \, \varepsilon^2.\] 
This implies 
\begin{align*} 
\bigg \| 
&\Omega \wedge F_A - \Omega_0 \wedge F_A + \sum_{k=1}^8 i_{e_k} \alpha(y) \wedge i_{e_k}(\Omega
_0 \wedge F_A) \\ 
&- \Omega_0 \wedge \sum_{k=1}^8 i_{e_k} \alpha(y) \wedge i_{e_k} F_A \bigg \|_{\mathcal{C}_2
^\gamma(M)} \leq C \, \varepsilon^2, 
\end{align*} 
hence
\begin{align*} 
\bigg \| 
&*(\Omega \wedge F_A) - *(\Omega_0 \wedge F_A) + \sum_{k=1}^8 i_{e_k} \alpha(y) \wedge i_{e_k} 
* (\Omega_0 \wedge F_A) \\ 
&- * \bigg ( \Omega_0 \wedge \sum_{k=1}^8 i_{e_k} \alpha(y) \wedge i_{e_k} F_A \bigg ) \bigg \|
_{\mathcal{C}_2^\gamma(M)} \leq C \, \varepsilon^2. 
\end{align*}  
Therefore, we obtain 
\begin{align*} 
\bigg \| 
&(F_A + *(\Omega \wedge F_A)) - (F_A + *(\Omega_0 \wedge F_A)) + \sum_{k=1}^8 i_{e_k} \alpha(y) 
\wedge i_{e_k} (F_A + * (\Omega_0 \wedge F_A)) \\ 
&- \sum_{k=1}^8 i_{e_k} \alpha(y) \wedge i_{e_k} F_A - * \bigg ( \Omega_0 \wedge \sum_{k=1}^8 i
_{e_k} \alpha(y) \wedge i_{e_k} F_A \bigg ) \bigg \|_{\mathcal{C}_2^\gamma(M)} \leq C \, 
\varepsilon^2. 
\end{align*} 
According to Proposition 3.4, we have 
\[\|F_A + * (\Omega_0 \wedge F_A)\|_{\mathcal{C}_3^\gamma(M)} \leq C \, \varepsilon^2,\] 
hence 
\[\bigg \| \sum_{k=1}^8 i_{e_k} \alpha(y) \wedge i_{e_k}(F_A + * (\Omega_0 \wedge F_A)) \bigg \|
_{\mathcal{C}_2^\gamma(M)} \leq C \, \varepsilon^2.\] 
Moreover, we have 
\[3 \, \sum_{k=1}^8 i_{e_k} \alpha(y) \wedge i_{e_k} F_A - * \bigg ( \Omega_0 \wedge \sum_{k=1}
^8 i_{e_k} \alpha(y) \wedge i_{e_k} F_A \bigg ) = 0.\] 
Thus, we conclude that 
\[\bigg \| (F_A + *(\Omega \wedge F_A)) - (F_A + *(\Omega_0 \wedge F_A)) - 4 \, \sum_{k=1}^8 i
_{e_k} \alpha(y) \wedge i_{e_k} F_A \bigg \|_{\mathcal{C}_2^\gamma(M)} \leq C \, \varepsilon^2.
\] 
The assertion follows now from Proposition 6.1. \\

\begin{proposition}
The curvature of $\tilde{A}$ satisfies 
\begin{align*} 
\bigg \| &\Pi(F_{\tilde{A}} + *(\Omega \wedge F_{\tilde{A}})) \\ 
&- 4 \, \bigg ( \proj_V \Big ( \sum_{i,j=1}^4 (\nabla_i v_k + \alpha_{ik,l} \, v_l) \, e_i 
\otimes e_k^\perp \Big ), \\ 
&\hspace{10.1mm} \proj_W \Big ( \sum_{i,k,l=1}^4 (\lambda^{-1} \, \nabla_i \lambda \, \delta
_{kl} + \theta_{i,kl} + \alpha_{ik,l}) \, e_i \otimes e_k^\perp \otimes e_l^\perp \bigg ) \bigg 
\|_{\mathcal{C}^\gamma(S)} \leq C \, \varepsilon^{\frac{1}{32}}. 
\end{align*}
\end{proposition}

\textbf{Proof.} 
Using the estimate 
\[\|a\|_{C_{2+\nu}^{1,\gamma}(M)} \leq C \, \varepsilon^{2-\nu-\gamma},\] 
we obtain 
\begin{align*} 
\big \| \Pi \big ( [a,a] + *(\Omega \wedge [a,a]) \big ) \big \|_{\mathcal{C}^\gamma(S)} 
&\leq C \, \varepsilon^{-2-\nu-\gamma} \, \|[a,a] + *(\Omega \wedge [a,a])\|_{\mathcal{C}_{3+
\nu}^\gamma(M)} \\ 
&\leq C \, \varepsilon^{-3-2\nu-\gamma} \, \|a\|_{\mathcal{C}_{2+\nu}^{1,\gamma}(M)}^2 \\ 
&\leq C \, \varepsilon^{1-4\nu-3\gamma}. 
\end{align*} 
Moreover, we have 
\[\|\Pi \, \mathbb{L}_A a\|_{\mathcal{C}^\gamma(S)} \leq C \, \varepsilon^{\frac{1}{32}}.\] 
Hence, the assertion follows from Proposition 6.2. \\

\textbf{Proof of Theorem 1.1.} 
Let 
\[\Xi_\varepsilon(v,\lambda,J,\omega) = \Pi(F_{\tilde{A}} + *(\Omega \wedge F_{\tilde{A}})).\] 
The first part of Theorem 1.1 follows from Proposition 5.2, the second part from Proposition 6.3. 
\\

\section{Discussion}

In this final section, we show how the first order balancing condition derived in this paper is 
related to the second order balancing condition in \cite{Br}. To this end, we assume that 
$\Omega$ is parallel. Then the Riemann curvature tensor of $M$ belongs to $\Lambda_-^2 TM 
\otimes \Lambda_-^2 TM$. Since $S$ is a Cayley submanifold, the second fundamental form of $S$ 
satisfies 
\begin{align*} 
&h(e_k,e_1,e_1^\perp) + h(e_k,e_2,e_2^\perp) + h(e_k,e_3,e_3^\perp) + h(e_k,e_4,e_4^\perp) = 0 
\\ 
&h(e_k,e_1,e_2^\perp) - h(e_k,e_2,e_1^\perp) - h(e_k,e_3,e_4^\perp) + h(e_k,e_4,e_3^\perp) = 0 
\\ 
&h(e_k,e_1,e_3^\perp) + h(e_k,e_2,e_4^\perp) - h(e_k,e_3,e_1^\perp) - h(e_k,e_4,e_2^\perp) = 0 
\\ 
&h(e_k,e_1,e_4^\perp) - h(e_k,e_2,e_3^\perp) + h(e_k,e_3,e_2^\perp) - h(e_k,e_4,e_1^\perp) = 0. 
\end{align*} 
We denote the curvature of the normal bundle $NS$ by $E$. Using the Gauss equations, we obtain 
\begin{align*} 
E(e_i,e_j,e_k^\perp,e_l^\perp) 
&= R(e_i,e_j,e_k^\perp,e_l^\perp) \\ 
&- \sum_{m=1}^4 h(e_m,e_i,e_k^\perp) \, h(e_m,e_j,e_l^\perp) + h(e_m,e_i,e_l^\perp) \, h(e_m,e
_j,e_k^\perp). 
\end{align*} 

Since $\nabla \Omega = 0$, the first part of the balancing condition becomes 
\begin{align*} 
&\nabla_1 v_1 + \nabla_2 v_2 + \nabla_3 v_3 + \nabla_4 v_4 = 0 \\ 
&\nabla_1 v_2 - \nabla_2 v_1 - \nabla_3 v_4 + \nabla_4 v_3 = 0 \\ 
&\nabla_1 v_3 + \nabla_2 v_4 - \nabla_3 v_1 - \nabla_4 v_2 = 0 \\ 
&\nabla_1 v_4 - \nabla_2 v_3 + \nabla_3 v_2 - \nabla_4 v_1 = 0. 
\end{align*} 
This implies 
\begin{align*} 
0 &= \Delta v_1 \\ 
&+ \nabla_1 \nabla_2 v_2 - \nabla_2 \nabla_1 v_2 + \nabla_3 \nabla_4 v_2 - \nabla_4 \nabla_3 v_2 
\\ 
&+ \nabla_1 \nabla_3 v_3 - \nabla_3 \nabla_1 v_3 + \nabla_4 \nabla_2 v_3 - \nabla_2 \nabla_4 v_3 
\\ 
&+ \nabla_1 \nabla_4 v_4 - \nabla_4 \nabla_1 v_4 + \nabla_2 \nabla_3 v_4 - \nabla_3 \nabla_2 v_4 
\end{align*} 
\begin{align*} 
0 &= \Delta v_2 \\ 
&- \nabla_1 \nabla_2 v_1 + \nabla_2 \nabla_1 v_1 - \nabla_3 \nabla_4 v_1 + \nabla_4 \nabla_3 v_1 
\\ 
&- \nabla_1 \nabla_3 v_4 + \nabla_3 \nabla_1 v_4 - \nabla_4 \nabla_2 v_4 + \nabla_2 \nabla_4 v_4 
\\ 
&+ \nabla_1 \nabla_4 v_3 - \nabla_4 \nabla_1 v_3 + \nabla_2 \nabla_3 v_3 - \nabla_3 \nabla_2 v_3 
\end{align*} 
\begin{align*} 
0 &= \Delta v_3 \\ 
&+ \nabla_1 \nabla_2 v_4 - \nabla_2 \nabla_1 v_4 + \nabla_3 \nabla_4 v_4 - \nabla_4 \nabla_3 v_4 
\\ 
&- \nabla_1 \nabla_3 v_1 + \nabla_3 \nabla_1 v_1 - \nabla_4 \nabla_2 v_1 + \nabla_2 \nabla_4 v_1 
\\ 
&- \nabla_1 \nabla_4 v_2 + \nabla_4 \nabla_1 v_2 - \nabla_2 \nabla_3 v_2 + \nabla_3 \nabla_2 v_2 
\end{align*} 
\begin{align*} 
0 &= \Delta v_4 \\ 
&- \nabla_1 \nabla_2 v_3 + \nabla_2 \nabla_1 v_3 - \nabla_3 \nabla_4 v_3 + \nabla_4 \nabla_3 v_3 
\\ 
&+ \nabla_1 \nabla_3 v_2 - \nabla_3 \nabla_1 v_2 + \nabla_4 \nabla_2 v_2 - \nabla_2 \nabla_4 v_2 
\\ 
&- \nabla_1 \nabla_4 v_1 + \nabla_4 \nabla_1 v_1 - \nabla_2 \nabla_3 v_1 + \nabla_3 \nabla_2 v_1. 
\end{align*} 
From this it follows that 
\begin{align*} 
0 &= \Delta v_1 \\ 
&+ \big ( E(e_1,e_2,e_1^\perp,e_2^\perp) + E(e_3,e_4,e_1^\perp,e_2^\perp) + E(e_1,e_3,e_1^\perp,
e_3^\perp) + E(e_4,e_2,e_1^\perp,e_3^\perp) \\ 
&+ E(e_1,e_4,e_1^\perp,e_4^\perp) + E(e_2,e_3,e_1^\perp,e_4^\perp) \big ) \, v_1 \\ 
&+ \big ( E(e_1,e_3,e_2^\perp,e_3^\perp) + E(e_4,e_2,e_2^\perp,e_3^\perp) + E(e_1,e_4,e_2^\perp,
e_4^\perp) + E(e_2,e_3,e_2^\perp,e_4^\perp) \big ) \, v_2 \\ 
&+ \big ( E(e_1,e_2,e_3^\perp,e_2^\perp) + E(e_3,e_4,e_3^\perp,e_2^\perp) + E(e_1,e_4,e_3^\perp,
e_4^\perp) + E(e_2,e_3,e_3^\perp,e_4^\perp) \big ) \, v_3 \\ 
&+ \big ( E(e_1,e_2,e_4^\perp,e_2^\perp) + E(e_3,e_4,e_4^\perp,e_2^\perp) + E(e_1,e_3,e_4^\perp,
e_3^\perp) + E(e_4,e_2,e_4^\perp,e_3^\perp) \big ) \, v_4 
\end{align*} 
\begin{align*} 
0 &= \Delta v_2 \\ 
&+ \big ( E(e_3,e_1,e_1^\perp,e_4^\perp) + E(e_2,e_4,e_1^\perp,e_4^\perp) + E(e_1,e_4,e_1^\perp,
e_3^\perp) + E(e_2,e_3,e_1^\perp,e_3^\perp) \big ) \, v_1 \\ 
&+ \big ( E(e_2,e_1,e_2^\perp,e_1^\perp) + E(e_4,e_3,e_2^\perp,e_1^\perp) + E(e_3,e_1,e_2^\perp,
e_4^\perp) + E(e_2,e_4,e_2^\perp,e_4^\perp) \\ 
&+ E(e_1,e_4,e_2^\perp,e_3^\perp) + E(e_2,e_3,e_2^\perp,e_3^\perp) \big ) \, v_2 \\ 
&+ \big ( E(e_2,e_1,e_3^\perp,e_1^\perp) + E(e_4,e_3,e_3^\perp,e_1^\perp) + E(e_3,e_1,e_3^\perp,
e_4^\perp) + E(e_2,e_4,e_3^\perp,e_4^\perp) \big ) \, v_3 \\ 
&+ \big ( E(e_2,e_1,e_4^\perp,e_1^\perp) + E(e_4,e_3,e_4^\perp,e_1^\perp) + E(e_1,e_4,e_4^\perp,
e_3^\perp) + E(e_2,e_3,e_4^\perp,e_3^\perp) \big ) \, v_4 
\end{align*} 
\begin{align*} 
0 &= \Delta v_3 \\ 
&+ \big ( E(e_1,e_2,e_1^\perp,e_4^\perp) + E(e_3,e_4,e_1^\perp,e_4^\perp) + E(e_4,e_1,e_1^\perp,
e_2^\perp) + E(e_3,e_2,e_1^\perp,e_2^\perp) \big ) \, v_1 \\ 
&+ \big ( E(e_1,e_2,e_2^\perp,e_4^\perp) + E(e_3,e_4,e_2^\perp,e_4^\perp) + E(e_3,e_1,e_2^\perp,
e_1^\perp) + E(e_2,e_4,e_2^\perp,e_1^\perp) \big ) \, v_2 \\ 
&+ \big ( E(e_1,e_2,e_3^\perp,e_4^\perp) + E(e_3,e_4,e_3^\perp,e_4^\perp) + E(e_3,e_1,e_3^\perp,
e_1^\perp) + E(e_2,e_4,e_3^\perp,e_1^\perp) \\ 
&+ E(e_4,e_1,e_3^\perp,e_2^\perp) + E(e_3,e_2,e_3^\perp,e_2^\perp) \big ) \, v_3 \\ 
&+ \big ( E(e_3,e_1,e_4^\perp,e_1^\perp) + E(e_2,e_4,e_4^\perp,e_1^\perp) + E(e_4,e_1,e_4^\perp,
e_2^\perp) + E(e_3,e_2,e_4^\perp,e_2^\perp) \big ) \, v_4
\end{align*} 
\begin{align*} 
0 &= \Delta v_4 \\ 
&+ \big ( E(e_2,e_1,e_1^\perp,e_3^\perp) + E(e_4,e_3,e_1^\perp,e_3^\perp) + E(e_1,e_3,e_1^\perp,
e_2^\perp) + E(e_4,e_2,e_1^\perp,e_2^\perp) \big ) \, v_1 \\ 
&+ \big ( E(e_2,e_1,e_2^\perp,e_3^\perp) + E(e_4,e_3,e_2^\perp,e_3^\perp) + E(e_4,e_1,e_2^\perp,
e_1^\perp) + E(e_3,e_2,e_2^\perp,e_1^\perp) \big ) \, v_2 \\ 
&+ \big ( E(e_1,e_3,e_3^\perp,e_2^\perp) + E(e_4,e_2,e_3^\perp,e_2^\perp) + E(e_4,e_1,e_3^\perp,
e_1^\perp) + E(e_3,e_2,e_3^\perp,e_1^\perp) \big ) \, v_3 \\ 
&+ \big ( E(e_2,e_1,e_4^\perp,e_3^\perp) + E(e_4,e_3,e_4^\perp,e_3^\perp) + E(e_1,e_3,e_4^\perp,
e_2^\perp) + E(e_4,e_2,e_4^\perp,e_2^\perp) \\ 
&+ E(e_4,e_1,e_4^\perp,e_1^\perp) + E(e_3,e_2,e_4^\perp,e_1^\perp) \big ) \, v_4. 
\end{align*} 
Hence, we obtain 
\[0 = \Delta v_1 + \sum_{i,j,k=1}^4 h(e_i,e_j,e_1^\perp) \, h(e_i,e_j,e_k^\perp) \, v_k + \sum
_{i,k=1}^4 R(e_i,e_1^\perp,e_k^\perp,e_i) \, v_k\] 
\[0 = \Delta v_2 + \sum_{i,j,k=1}^4 h(e_i,e_j,e_2^\perp) \, h(e_i,e_j,e_k^\perp) \, v_k + \sum
_{i,k=1}^4 R(e_i,e_2^\perp,e_k^\perp,e_i) \, v_k\] 
\[0 = \Delta v_3 + \sum_{i,j,k=1}^4 h(e_i,e_j,e_3^\perp) \, h(e_i,e_j,e_k^\perp) \, v_k + \sum
_{i,k=1}^4 R(e_i,e_3^\perp,e_k^\perp,e_i) \, v_k\] 
\[0 = \Delta v_4 + \sum_{i,j,k=1}^4 h(e_i,e_j,e_4^\perp) \, h(e_i,e_j,e_k^\perp) \, v_k + \sum
_{i,k=1}^4 R(e_i,e_4^\perp,e_k^\perp,e_i) \, v_k.\] 

Furthermore, the second part of the balancing condition can be written in the form 
\begin{align*} 
&2\lambda^{-1} \, \nabla_1 \lambda + (\theta_{2,21} + \theta_{2,43}) + (\theta_{3,31} + \theta
_{3,24}) + (\theta_{4,41} + \theta_{4,32}) = 0 \\ 
&(\theta_{1,21} + \theta_{1,43}) - 2\lambda^{-1} \, \nabla_2 \lambda - (\theta_{3,41} + \theta
_{3,32}) + (\theta_{4,31} + \theta_{4,24}) = 0 \\ 
&(\theta_{1,31} + \theta_{1,24}) + (\theta_{2,41} + \theta_{2,32}) - 2\lambda^{-1} \, \nabla_3 
\lambda - (\theta_{4,21} + \theta_{4,43}) = 0 \\ 
&(\theta_{1,41} + \theta_{1,32}) - (\theta_{2,31} + \theta_{2,24}) + (\theta_{3,21} + \theta
_{3,43}) - 2\lambda^{-1} \, \nabla_4 \lambda = 0. 
\end{align*} 
This implies 
\begin{align*} 
0 &= 2\lambda^{-1} \, \Delta \lambda - 2\lambda^{-2} \, |\nabla \lambda|^2 \\ 
&+ \nabla_1 \theta_{2,21} - \nabla_2 \theta_{1,21} + \nabla_1 \theta_{2,43} - \nabla_2 \theta
_{1,43} \\ 
&+ \nabla_3 \theta_{4,21} - \nabla_4 \theta_{3,21} + \nabla_3 \theta_{4,43} - \nabla_4 \theta
_{3,43} \\ 
&+ \nabla_1 \theta_{3,31} - \nabla_3 \theta_{1,31} + \nabla_1 \theta_{3,24} - \nabla_3 \theta
_{1,24} \\ 
&+ \nabla_4 \theta_{2,31} - \nabla_2 \theta_{4,31} + \nabla_4 \theta_{2,24} - \nabla_2 \theta
_{4,24} \\ 
&+ \nabla_1 \theta_{4,41} - \nabla_4 \theta_{1,41} + \nabla_1 \theta_{4,32} - \nabla_4 \theta
_{1,32} \\ 
&+ \nabla_2 \theta_{3,41} - \nabla_3 \theta_{2,41} + \nabla_2 \theta_{3,32} - \nabla_3 \theta
_{2,32}, 
\end{align*} 
hence 
\begin{align*} 
0 &= 2\lambda^{-1} \, \Delta \lambda - \frac{1}{2} \, |\theta|^2 \\ 
&+ \nabla_1 \theta_{2,21} - \nabla_2 \theta_{1,21} + [\theta_1,\theta_2]_{21} + \nabla_1 \theta
_{2,43} - \nabla_2 \theta_{1,43} + [\theta_1,\theta_2]_{43} \\ 
&+ \nabla_3 \theta_{4,21} - \nabla_4 \theta_{3,21} + [\theta_3,\theta_4]_{21} + \nabla_3 \theta
_{4,43} - \nabla_4 \theta_{3,43} + [\theta_3,\theta_4]_{43} \\ 
&+ \nabla_1 \theta_{3,31} - \nabla_3 \theta_{1,31} + [\theta_1,\theta_3]_{31} + \nabla_1 \theta
_{3,24} - \nabla_3 \theta_{1,24} + [\theta_1,\theta_3]_{24} \\ 
&+ \nabla_4 \theta_{2,31} - \nabla_2 \theta_{4,31} + [\theta_4,\theta_2]_{31} + \nabla_4 \theta
_{2,24} - \nabla_2 \theta_{4,24} + [\theta_4,\theta_2]_{24} \\ 
&+ \nabla_1 \theta_{4,41} - \nabla_4 \theta_{1,41} + [\theta_1,\theta_4]_{41} + \nabla_1 \theta
_{4,32} - \nabla_4 \theta_{1,32} + [\theta_1,\theta_4]_{32} \\ 
&+ \nabla_2 \theta_{3,41} - \nabla_3 \theta_{2,41} + [\theta_2,\theta_3]_{41} + \nabla_2 \theta
_{3,32} - \nabla_3 \theta_{2,32} + [\theta_2,\theta_3]_{32}. 
\end{align*} 
From this it follows that 
\begin{align*} 
0 &= 2\lambda^{-1} \, \Delta \lambda - \frac{1}{2} \, |\theta|^2 \\ 
&+ E(e_1,e_2,e_1^\perp,e_2^\perp) + E(e_1,e_2,e_3^\perp,e_4^\perp) + E(e_3,e_4,e_1^\perp,e_2
^\perp) + E(e_3,e_4,e_3^\perp,e_4^\perp) \\ 
&+ E(e_1,e_3,e_1^\perp,e_3^\perp) + E(e_1,e_3,e_4^\perp,e_2^\perp) + E(e_4,e_2,e_1^\perp,e_3
^\perp) + E(e_4,e_2,e_4^\perp,e_2^\perp) \\ 
&+ E(e_1,e_4,e_1^\perp,e_4^\perp) + E(e_1,e_4,e_2^\perp,e_3^\perp) + E(e_2,e_3,e_1^\perp,e_4
^\perp) + E(e_2,e_3,e_2^\perp,e_3^\perp). 
\end{align*} 
Using the identities 
\begin{align*} 
&E(e_1,e_2,e_1^\perp,e_2^\perp) + E(e_3,e_4,e_1^\perp,e_2^\perp) + E(e_1,e_3,e_1^\perp,e_3
^\perp) + E(e_4,e_2,e_1^\perp,e_3^\perp) \\ 
&+ E(e_1,e_4,e_1^\perp,e_4^\perp) + E(e_2,e_3,e_1^\perp,e_4^\perp) \\ 
&= \sum_{i,j=1}^4 h(e_i,e_j,e_1^\perp) \, h(e_i,e_j,e_1^\perp) + \sum_{i=1}^4 R(e_i,e_1^\perp,e
_1^\perp,e_i) 
\end{align*}
\begin{align*} 
&E(e_2,e_1,e_2^\perp,e_1^\perp) + E(e_4,e_3,e_2^\perp,e_1^\perp) + E(e_3,e_1,e_2^\perp,e_4
^\perp) + E(e_2,e_4,e_2^\perp,e_4^\perp) \\ 
&+ E(e_1,e_4,e_2^\perp,e_3^\perp) + E(e_2,e_3,e_2^\perp,e_3^\perp) \\ 
&= \sum_{i,j=1}^4 h(e_i,e_j,e_2^\perp) \, h(e_i,e_j,e_2^\perp) + \sum_{i=1}^4 R(e_i,e_2^\perp,e
_2^\perp,e_i) 
\end{align*} 
\begin{align*} 
&E(e_1,e_2,e_3^\perp,e_4^\perp) + E(e_3,e_4,e_3^\perp,e_4^\perp) + E(e_3,e_1,e_3^\perp,e_1
^\perp) + E(e_2,e_4,e_3^\perp,e_1^\perp) \\ 
&+ E(e_4,e_1,e_3^\perp,e_2^\perp) + E(e_3,e_2,e_3^\perp,e_2^\perp) \\ 
&= \sum_{i,j=1}^4 h(e_i,e_j,e_3^\perp) \, h(e_i,e_j,e_3^\perp) + \sum_{i=1}^4 R(e_i,e_3^\perp,e
_3^\perp,e_i) 
\end{align*} 
\begin{align*} 
&E(e_2,e_1,e_4^\perp,e_3^\perp) + E(e_4,e_3,e_4^\perp,e_3^\perp) + E(e_1,e_3,e_4^\perp,e_2
^\perp) + E(e_4,e_2,e_4^\perp,e_2^\perp) \\ 
&+ E(e_4,e_1,e_4^\perp,e_1^\perp) + E(e_3,e_2,e_4^\perp,e_1^\perp) \\ 
&= \sum_{i,j=1}^4 h(e_i,e_j,e_4^\perp) \, h(e_i,e_j,e_4^\perp) + \sum_{i=1}^4 R(e_i,e_4^\perp,e
_4^\perp,e_i), 
\end{align*} 
we obtain 
\begin{align*} 
&2 \, E(e_1,e_2,e_1^\perp,e_2^\perp) + 2 \, E(e_1,e_2,e_3^\perp,e_4^\perp) + 2 \, E(e_3,e_4,e_1
^\perp,e_2^\perp) + 2 \, E(e_3,e_4,e_3^\perp,e_4^\perp) \\ 
&+ 2 \, E(e_1,e_3,e_1^\perp,e_3^\perp) + 2 \, E(e_1,e_3,e_4^\perp,e_2^\perp) + 2 \, E(e_4,e_2,e
_1^\perp,e_3^\perp) + 2 \, E(e_4,e_2,e_4^\perp,e_2^\perp) \\ 
&+ 2 \, E(e_1,e_4,e_1^\perp,e_4^\perp) + 2 \, E(e_1,e_4,e_2^\perp,e_3^\perp) + 2 \, E(e_2,e_3,e
_1^\perp,e_4^\perp) + 2 \, E(e_2,e_3,e_2^\perp,e_3^\perp) \\ 
&= \sum_{i,j,k=1}^4 h(e_i,e_j,e_k^\perp) \, h(e_i,e_j,e_k^\perp) + \sum_{i,k=1}^4 R(e_i,e_k
^\perp,e_k^\perp,e_i). 
\end{align*} 
Thus, we conclude that 
\[0 = 4\lambda^{-1} \, \Delta \lambda - |\theta|^2 + \sum_{i,j,k=1}^4 h(e_i,e_j,e_k^\perp) \, 
h(e_i,e_j,e_k^\perp) + \sum_{i,k=1}^4 R(e_i,e_k^\perp,e_k^\perp,e_i).\] 
A similar calculation gives 
\[0 = \sum_{i=1}^4 (\nabla_i \theta_i + 2\lambda^{-1} \, \nabla_i \lambda \, \theta_i).\] 
Thus, the first order balancing condition implies the second order balancing condition derived 
in \cite{Br}. \\


\begin{thebibliography}{99}
\bibitem{Ba}
A. Bahri, \textit{An invariant for Yamabe-type flows with applications to scalar curvature 
problems in high dimension,} Duke Math. J. 81, 323-466 (1996)

\bibitem{BC}
A. Bahri and J.M. Coron, \textit{On a nonlinear elliptic equation involving the critical Sobolev 
exponent: The effect of the topology of the domain,} Comm. Pure Appl. Math. 41, 253-290 (1988)

\bibitem{BKS}
L. Baulieu, H. Kanno, and I. M. Singer, \textit{Special quantum field theories in eight and 
other dimensions,} Comm. Math. Phys. 194, 149-175 (1998)

\bibitem{BBMOOY}
K. Becker, M. Becker, D. Morrison, H. Ooguri, Y. Oz, and Z. Yin, \textit{Supersymmetric cycles 
in exceptional holonomy manifolds and Calabi-Yau $4$-folds,} Nucl. Phys. B 480, 225-238 (1996)

\bibitem{Br}
S. Brendle, \textit{On the construction of solutions to the Yang-Mills equations in higher 
dimensions,} preprint (2002)

\bibitem{Ch}
J. Chen, \textit{Complex anti-self-dual connections on a product of Calabi-Yau surfaces and 
triholomorphic curves,} Comm. Math. Phys. 201, 217-247 (1999)

\bibitem{DK}
S. K. Donaldson and P. B. Kronheimer, \textit{The Geometry of Four-Manifolds,} Oxford University 
Press (1990)

\bibitem{DT}
S. K. Donaldson and R. P. Thomas, \textit{Gauge theory in higher dimensions,} The geometric 
universe (Oxford 1996), 31-47, Oxford University Press (1998)

\bibitem{HL}
R. Harvey and H. B. Lawson, \textit{Calibrated Geometries,} Acta Math. 148, 47-157 (1982)

\bibitem{Hi}
N. Hitchin, \textit{Lectures on special Lagrangian submanifolds,} Winter School on Mirror 
Symmetry, Vector Bundles and Lagrangian Submanifolds (Cambridge, MA, 1999), 151-182, AMS/IP 
Stud. Adv. Math. 23, Amer. Math. Soc. Providence, RI, 2001

\bibitem{Jo}
D. Joyce, \textit{Compact $8$-manifolds with holonomy $Spin(7)$,} Invent. Math. 123, 507-552 
(1996)

\bibitem{Ka1}
N. Kapouleas, \textit{Complete constant mean curvature surfaces in Euclidean three-space,} Ann. 
of Math. 131, 239-330 (1990)

\bibitem{Ka2}
N. Kapouleas, \textit{Compact constant mean curvature surfaces in Euclidean three-space,} J. 
Diff. Geom. 33, 683-715 (1991)

\bibitem{KMP}
R. Kusner, R. Mazzeo, and D. Pollack, \textit{The moduli space of of complete embedded constant 
mean curvature surfaces,} Geom. Funct. Anal. 6, 120-137 (1996)

\bibitem{Li}
F.-H. Lin, \textit{Complex Ginzburg-Landau equations and dynamics of vortices, filaments, and 
codimension-$2$ submanifolds,} Comm. Pure Appl. Math. 51, 385-441 (1998)

\bibitem{LR}
F.-H. Lin and T. Rivi\`ere, \textit{Complex Ginzburg-Landau equations in high dimensions and 
codimension two area minimizing currents,} J. Eur. Math. Soc. 1, 237-311 (1999)

\bibitem{MP1}
R. Mazzeo and F. Pacard, \textit{A construction of singular solutions for a semilinear elliptic 
equation using asymptotic analysis,} J. Diff. Geom. 44, 331-370 (1996)

\bibitem{MP2}
R. Mazzeo and F. Pacard, \textit{Constant mean curvature surfaces with Delaunay ends,} Comm. 
Anal. Geom. 9, 169-237 (2001)

\bibitem{MPP}
R. Mazzeo, F. Pacard and D. Pollack, \textit{Connected sums of constant mean curvature surfaces 
in Euclidean $3$-space,} J. Reine Angew. Math. 536, 115-165 (2001)

\bibitem{MS}
R. Mazzeo and N. Smale, \textit{Conformally flat metrics of constant positive scalar curvature 
on subdomains of the sphere,} J. Diff. Geom. 34, 581-621 (1991)

\bibitem{PR1}
F. Pacard and M. Ritor\'e, \textit{From constant mean curvature hypersurfaces to the gradient 
theory of phase transitions,} preprint (2003)

\bibitem{PR2}
F. Pacard and T. Rivi\`ere, \textit{Linear and Nonlinear Aspects of Vortices. The 
Ginzburg-Landau Model,} Progress in Nonlinear Differential Equations and their Applications, 
vol. 39, Birkh\"auser, Boston (2000)

\bibitem{Ta1}
C. H. Taubes, \textit{Self-dual Yang-Mills connections on non-self-dual $4$-manifolds,} J. Diff. 
Geom. 17, 139-170 (1982)

\bibitem{Ta2}
C. H. Taubes, \textit{$\text{SW} \Longrightarrow \text{Gr}$: from the Seiberg-Witten equations 
to pseudo-holomorphic curves,} J. Amer. Math. Soc. 9, 845-918 (1996)

\bibitem{Ta3}
C. H. Taubes, \textit{$\text{Gr} \Longrightarrow \text{SW}$: from pseudo-holomorphic curves to 
Seiberg-Witten solutions,} J. Diff. Geom. 51, 203-334 (1999)

\bibitem{Th}
R. P. Thomas, \textit{A holomorphic Casson invariant for Calabi-Yau $3$-folds and bundles on 
$K3$ fibrations,} J. Diff. Geom. 54, 367-438 (2000)

\bibitem{Ti}
G. Tian, \textit{Gauge theory and calibrated geometry,} Ann. of Math. 151, 193-268 (2000)
\end{thebibliography}
\end{document}